\RequirePackage{fix-cm}
\documentclass[smallcondensed]{svjour3}     
\smartqed  
\journalname{Journal of Scientific Computing}
\usepackage{amsopn}	

\usepackage{amssymb}
\usepackage{amsmath}

\usepackage{amsthm}
\usepackage{subcaption}
\usepackage{color}
\usepackage{marginnote}
\usepackage{algorithm} 
\usepackage{algorithmic}
\usepackage{multirow}

\usepackage{graphicx}

\graphicspath{ {./images/} }

\newcommand{\E}[2][]{\mathbb{E}_{{#1}} \left[ {#2} \right] }
\newcommand{\EMC}[2][]{\mathbb{E_{MC}}_{{#1}} \left[ {#2} \right] }
\newcommand{\EML}[2][]{\mathbb{E_{ML}}_{{#1}} \left[ {#2} \right] }

\newcommand{\dEML}[2][]{\big( \mathbb{E_{ML}}_{{#1}} - \mathbb{E}_{{#1}} \big) \left[ {#2} \right] }
\newcommand{\dEMC}[2][]{\big( \mathbb{E_{MC}}_{{#1}} - \mathbb{E}_{{#1}} \big) \left[ {#2} \right] }

\newcommand{\norm}[2]{ \left\| #1 \right\|_{#2} }
\newcommand{\change}[1]{#1}

\newtheorem{Theorem}{Theorem}

\everymath{\displaystyle}

\begin{document}

\title{Method of Green's potentials for elliptic PDEs in domains with random \change{apertures} \thanks{This manuscript has been authored by UT-Battelle, LLC, under contract DE-AC05-00OR22725 with the US Department of Energy (DOE). The US government retains and the publisher, by accepting the article for publication, acknowledges that the US government retains a nonexclusive, paid-up, irrevocable, worldwide license to publish or reproduce the published form of this manuscript, or allow others to do so, for US government purposes. DOE will provide public access to these results of federally sponsored research in accordance with the DOE Public Access Plan (http://energy.gov/downloads/doe-public-access-plan)} }



\author{Viktor Reshniak         \and
        Yuri Melnikov 
}


\institute{
    V. Reshniak \at
    Computer Science and Mathematics Division, Oak Ridge National Laboratory, Oak Ridge, TN 37831, USA \\
    \email{reshniakv@ornl.gov}           
    \and
    Yu. A. Melnikov \at
    Department of Mathematical Sciences, Middle Tennessee State University, Murfreesboro, TN 37132, USA \\
    \email{Yuri.Melnikov@mtsu.edu}  
}

\date{Received: date / Accepted: date}

\maketitle

\begin{abstract}
Problems with topological uncertainties appear in many fields ranging from nano-device engineering to the design of bridges. 
In many of such problems, a part of the domains boundaries is subjected to random perturbations making inefficient conventional schemes that rely on discretization of the whole domain. 
In this paper, we study elliptic PDEs in domains with boundaries comprised of \change{a deterministic part and random apertures}, and apply the method of modified potentials with Green's kernels defined on the deterministic part of the domain. 
This approach allows to reduce the dimension of the original differential problem by reformulating it as a boundary integral equation posed on \change{the random apertures} only. 
The multilevel Monte Carlo method is then applied to this modified integral equation \change{and its optimal $\epsilon^{-2}$ asymptotical complexity is shown}. 
\change{Finally,} we provide the qualitative analysis of the proposed technique and support it with numerical results.

\end{abstract}


\section{Introduction}
\label{intro}

Scientific and technological development is often linked to the increase in the requirements to the accuracy of mathematical models and numerical methods.
As an example, consideration of uncertainties in model inputs and parameters has been attracting a lot of attention of a research community in recent years and a large database of methods for the numerical treatment of boundary-value problems with random coefficients and random input data has been collected.
Problems with topological uncertainties have been also studied and their importance has been recognized in many applied fields ranging from nano-device engineering and analysis of micro electromechanical systems (MEMS) \cite{Agarwal2007,Zhu2004} to the design of bridges~\cite{Canuto2009}.
Other applications include flows over rough surfaces \cite{Tartakovsky2006a,Zayernouri2013}, surface imaging \cite{Tsong2005}, corrosion or wear of surfaces, homogenization of random heterogeneous media \cite{Savvas2014} and even modelling of blood flow \cite{Park2012}.

Existing numerical methods for PDEs in random domains differ by the way of approximating the spatial and random components of the solution.
For example, the random solution of a boundary value problem under the small noise assumption can be often represented as a sum of a deterministic component corresponding to the fixed nominal boundary and a small random perturbation which can be quantified using the methods of the ``shape calculus" \cite{Honda2005,Harbrecht2008}.

Alternatively, the original problem in a random domain can be transformed to the problem with random coefficients posed on a deterministic reference domain by means of the random change of variables estimated from a series of auxiliary PDEs \cite{Tartakovsky2006b}.
In conjunction with the stochastic Galerkin method, this approach was considered in \cite{Tartakovsky2006a,Kundu2014,Harbrecht2014} 
while the stochastic collocation approximation in random space was applied in \cite{Castrillon2014,Castrillon2015}. 
An equivalent Lagrangian approach was also proposed in \cite{Agarwal2007} where the mapping to the reference domain was combined with the stochastic spectral boundary element approximation. 

Similarly to the domain mapping method, the random displacement field can be applied directly to the mesh-based representation of the geometry producing the new mesh with random coordinates of nodes but the same fixed connectivity.
The main advantage of the mesh-based formulation is that the structure of the underlying linear system remains unchanged enabling reusability of the existing deterministic solvers.
This idea was proposed in \cite{Mohan2011} in combination with the polynomial chaos approximation in random space and later was studied in \cite{Harbrecht2014} in the context of the Quasi-Monte Carlo method applied to the random interface problem.

In fictitious domain methods, the original problem is reformulated on a larger deterministic domain containing all realizations of the random boundary.
The enclosing domain can be chosen arbitrarily allowing for simple discretizations which do not have to conform with the random boundaries at the cost of adding new variables to enforce the true boundary conditions.
For example, the authors of \cite{Nouy2007,Nouy2008} enriched the finite element approximation spaces with the suitably constructed functions which allow for the explicit representation of a solution in terms of the random variables describing the geometry.
In \cite{Savvas2014,Lang2013,Nouy2010}, this method was successfully applied to a number of problems with stochastic material interfaces.
The authors of \cite{Canuto2007} satisfied the boundary conditions by introducing the Lagrange multiplier which transformed the original elliptic equation into the larger saddle-point problem.
However, the information on the random geometry in the resulting linear system was encoded only in the part of the matrix coupling the primal variable and the boundary supported Lagrange multiplier.
Such localization property is a certain advantage of this approach over the domain mapping methods which propagate the boundary uncertainty to the whole domain.
Additionally, the method requires no assumptions on the size of the random displacements which is a major limitation of the perturbation techniques.

The aim of this paper is to construct an efficient and accurate numerical method with good localization properties in the sense described above.
Our motivation for such formulation is driven by the problems with only certain (often relatively small) part of the boundary being subjected to the random perturbations.
\change{In particular, the emphasis of this effort is on elliptic equations in arbitrary (otherwise deterministic) domains with random apertures.}
The fully discrete formulations of conventional solvers for such problems can lack efficiency due to the necessity in the discretization of the whole physical domain.
In this regard, the semi-analitical approximations hold a vast potential. 
Here we propose to adapt the method of Green's potentials for elliptic equations to the case of domains with random boundaries.
Green's potentials are the layer potentials with the modified kernels given by the suitably constructed Green's functions.
We will show that the proposed method allows to formulate the original boundary value problem in terms of the integral equations on the random part of the boundary only leading to the potentially significant computational savings.
The first step towards the practical application of this approach was done in \cite{Melnikov1977} where the so-called method of ``modified potentials" was introduced.
Later it was successfully applied to both stationary and time-dependent deterministic problems \cite{Melnikov2014,Melnikov1996}.
It is worth noting that the importance of Green's functions has been already recognized in various areas of the uncertainty quantification  \cite{Collins2015,BarajasSolano2013,Neuman1993,Neuman1996,Manolis1996}. 
Here we apply the Multilevel Monte Carlo (MLMC) method \cite{Giles2008} for the statistical approximation.
However, we note that any method of collocation type can be trivially adopted to the proposed numerical technique.

The paper is organized as follows.
In section \ref{sec:formulation}, we formulate the problem and introduce the equations and necessary analytical tools. 
In section \ref{sec:scheme}, we discuss the discretization scheme for the given equations.
The complexity analysis of the proposed scheme is given in section \ref{sec:complexity}.
Finally, the numerical examples in section \ref{sec:numerical} are provided in support of the obtained analytical results.

\section{Problem setting}
\label{sec:formulation}

Let $(\Omega,\mathcal{F},\mathbb{P})$ be a complete probability space with a set of outcomes $\Omega$, a sigma algebra of events $\mathcal{F}$ and a probability measure $\mathbb{P}$ defined on it.
For each outcome $\omega \in \Omega$, define $D(\omega)$ to be the realization of a random domain with a boundary comprised of deterministic and random parts $\partial D(\omega):=\partial D_1 \cup \partial D_2(\omega)$.
\change{$\partial D_1$ is assumed to be Lipschitz and $\partial D_2(\omega)$ is assumed to be analytic.}
We are concerned with the solutions of the following boundary value problem
\begin{alignat}{3} \label{eq:Poisson}
	\nonumber
	-\nabla^2 u(x,\omega)                                                    &= f(x)   &\qquad &\text{for } x \in D(\omega),
	\\ 
	\alpha_1 u(x,\omega) + \beta_1 \frac{\partial u(x,\omega)}{\partial n}   &= b_1(x)        &\qquad &\text{for } x \in \partial D_1,
	\\ \nonumber
	\alpha_2 u(x,\omega) + \beta_2 \frac{\partial u(x,\omega)}{\partial n}   &= b_2(x) &\qquad &\text{for } x \in \partial D_2(\omega).
\end{alignat}
We assume that, for each $\omega$, the solution $u(x,\omega)$ to the above problem exists, is unique and belongs to $H^1(D)$, the space of square integrable functions with square integrable first derivatives.
Additionally, we require $u({x},\omega)$ to be a Bochner integrable function with values in $H^1(D)$, i.e., $u({x}, \omega) \in L^p(\Omega; H^1(D))$, the function space given by
\begin{align*}
	L^p \Big( \Omega;&H^1(D) \Big) 
	\\ 
	&:=\left\{
	u : \Omega \to H^1(D) \; \Big| 
	\begin{array}{l}
		u \text{ is strongly measurable and } \norm{u}{{L^p(\Omega;H^1(D))}} < \infty	
	\end{array}		
	\right\}
\end{align*}
with the corresponding norm
\begin{align}\label{eq:rand_norm}
	\norm{u}{L^p(\Omega;H^1(D))}^p = 
	\left\{
		\begin{array}{l l}
			\displaystyle{ \int_{\Omega} \norm{ u(\cdot,\omega) }{H^1(D)}^p d \mathbb{P}(\omega) } &\qquad \text{if  } 0 < p < \infty,
			\\[1em]
			ess \sup_{\omega \in \Omega} \norm{ u(\cdot,\omega) }{H^1(D)} &\qquad \text{if  } p=\infty.
		\end{array}
	\right.
\end{align}
For compactness, we will use $L^p(\Omega)$ instead of $L^p(\Omega; H^1(D))$ in later discussions. 

\begin{figure}[t!]
	\centering
    \includegraphics[width=0.3\textwidth]{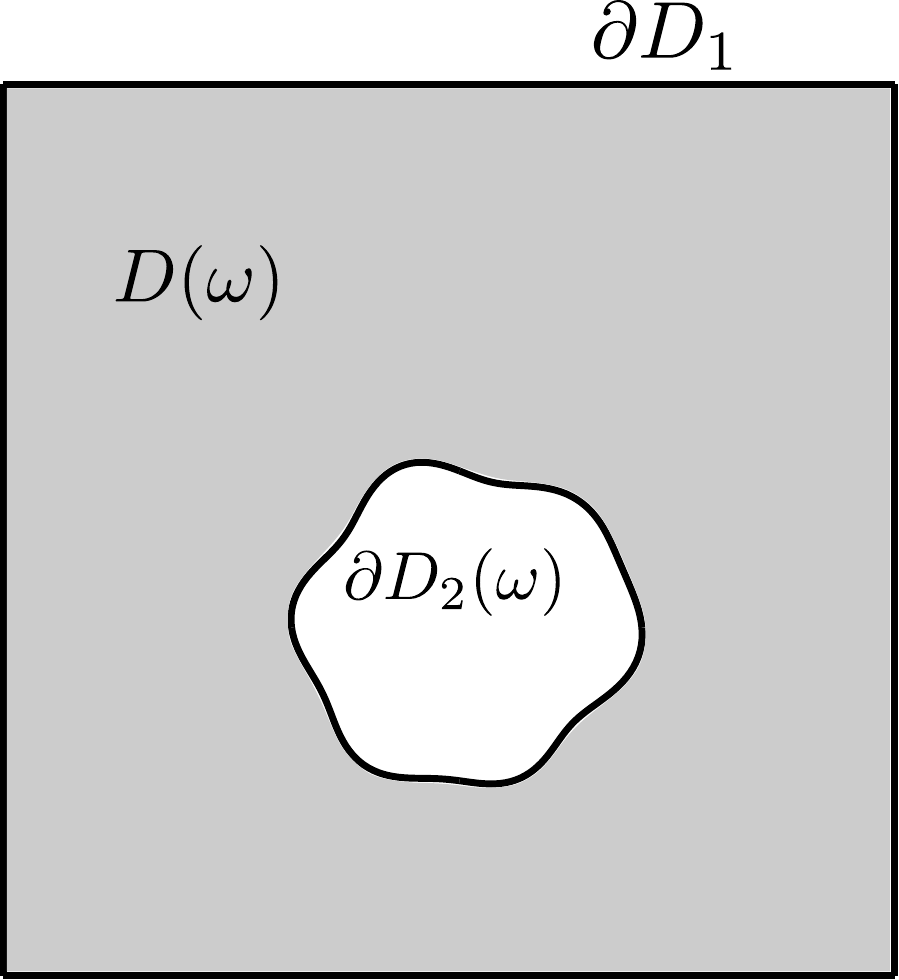}
    \qquad \hspace{5mm} \qquad
    \includegraphics[width=0.3\textwidth]{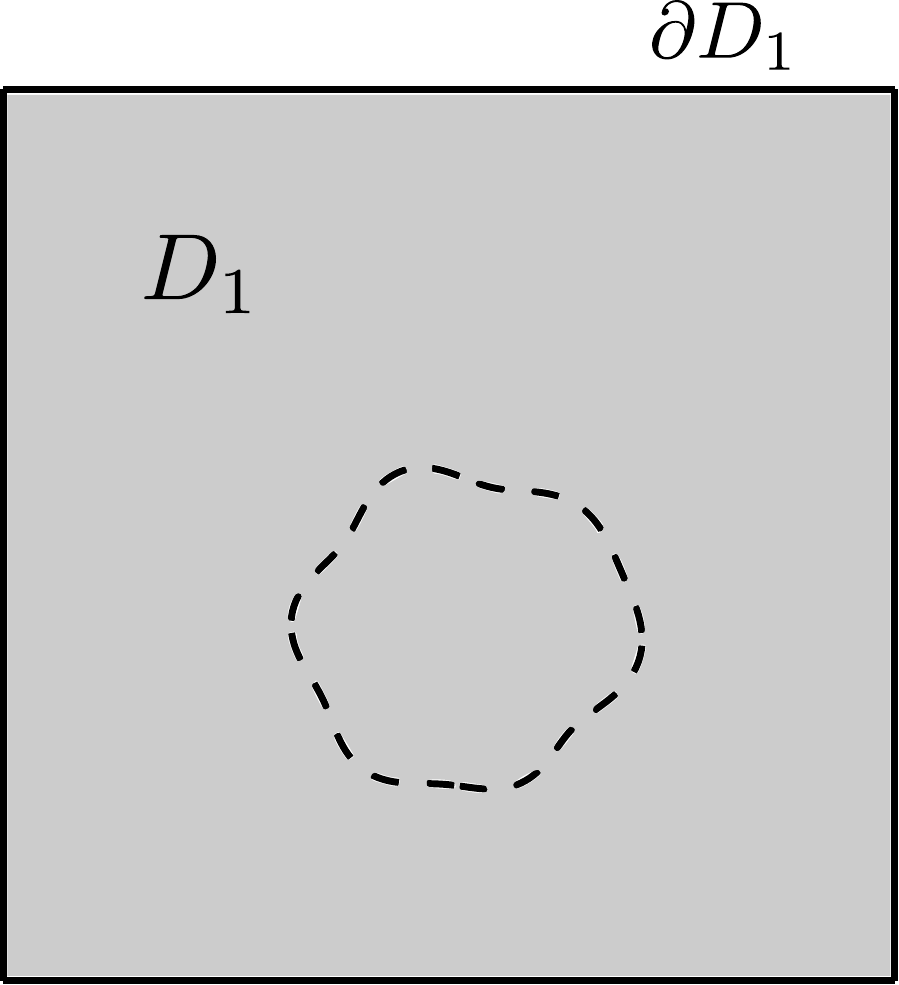}
	\caption[Realization of the random domain.]{ Realization of the random domain $D(\omega)$ (left) and the corresponding deterministic domain $D_1$ (right). }
	\label{fig:geom}
\end{figure}

Denote by $D_1$ the reference deterministic domain containing all realizations of the random boundary $\partial D_2(\omega)$.
This definition is similar to that used in the fictitious domain methods. 
However, we explicitly require that the deterministic part of the boundary $\partial D_1$ is also the boundary of $D_1$.
This definition is depicted graphically in Figure \ref{fig:geom}.
Due to linearity of the operators in \eqref{eq:Poisson}, one can represent the solution $u(x,\omega)$ as a superposition of two functions 
\begin{align}\label{eq:superpos}
	u(x,\omega) = u_1(x) + u_2(x,\omega),
\end{align}
where $u_1(x)$ is the deterministic component which satisfies the boundary value problem in the reference domain $D_1$
\begin{alignat}{3} \label{eq:u_1}
	-\nabla^2 u_1(x)                                               & = f(x)  &\qquad &\text{for } x \in D_1,
	\\ \nonumber
	\alpha_1 u_1(x) + \beta_1 \frac{\partial u_1(x)}{\partial n}   & = b_1(x)  &\qquad &\text{for } x \in \partial D_1
\end{alignat}
and the random component $u_2(x,\omega)$ can be determined from the following homogeneous boundary value problem
\begin{alignat}{3} \label{eq:u_2}
	\nonumber
	-\nabla^2 u_2(x,\omega)                                                      & = 0  &\qquad &\text{for } x \in D(\omega),
	\\ 
	\alpha_1 u_2(x,\omega) + \beta_1 \frac{\partial u_2(x,\omega)}{\partial n}   & = 0  &\qquad &\text{for } x \in \partial D_1,
	\\ \nonumber
	\alpha_2 u_2(x,\omega) + \beta_2 \frac{\partial u_2(x,\omega)}{\partial n}   & = \phi(x)  &\qquad &\text{for } x \in \partial D_2(\omega)
\end{alignat}
with the boundary condition on $\partial D_2(\omega)$ defined by the trace of $u_1(x)$ as follows
\begin{align} \label{eq:BC}
	\phi(x) = b_2(x) - \alpha_2 u_1(x) - \beta_2 \frac{\partial u_1(x)}{\partial n}.
\end{align}

\subsection{Method of Green's potentials}

The particularly simple structure of the problem in \eqref{eq:u_2} is very well suited for the construction of efficient solvers.
One of such solvers, namely the method of Green's potentials, is proposed in this section.

Firstly, define the Green's function corresponding to the boundary value problem \eqref{eq:u_1} as a solution of the complementary problem
\begin{alignat}{3} \label{eq:GF_2}
	-\nabla^2 G_1(x,\xi)                                                     & = \delta(\xi)  &\qquad &\text{for } x, \xi \in D_1,
	\\ \label{eq:GF_BC}
	\alpha_1 G_1(x,\xi) + \beta_1 \frac{\partial G_1(x,\xi)}{\partial n_x}   & = 0  &\qquad &\text{for } x \in \partial D_1,
\end{alignat}
where $\delta(\xi)$ is the Dirac measure of unit mass at point $\xi$.
With this function at hand, the solution $u_2(x,\omega)$ of the problem in \eqref{eq:u_2} allows the representation in the form of the single-layer potential \cite{Jaswon1963,Symm1963}
\begin{align} \label{eq:potential}
	\nonumber
	u_2(x,\omega) &= \int\displaylimits_{\partial D_1} G_1(x,y) \nu_1(y) d l(y) + \int\displaylimits_{\partial D_2(\omega)} G_1(x,y(\omega)) \nu_2(y(\omega)) d l(y(\omega))
	\\
	&= \int\displaylimits_{\partial D_2(\omega)} G_1(x,y(\omega)) \nu_2(y(\omega)) d l(y(\omega))
	= \int\displaylimits_0^1 G_1(x,\eta_{\omega}(t)) \mu(\eta_{\omega}(t)) dt,
\end{align}
where $\mu(\eta_{\omega}(t)) = \nu_2(\eta_{\omega}(t)) \left| \eta_{\omega}'(t) \right|$ and $\eta_{\omega}(t)$ defines a \change{smooth} parameterization of the random boundary curve.
Note that in the above expression, the first integral over the deterministic part of the boundary vanishes because of the Green's function $G_1(x,\xi)$.

Using \eqref{eq:BC} and the jump conditions of the derivative of the single-layer potential on the boundary, the unknown density $\mu(\eta_{\omega}(t))$ of the potential in \eqref{eq:potential} can be obtained from the following boundary integral equation
\begin{align} \label{eq:BIE}
    - \frac{\beta_2}{2} \frac{\mu(\eta_{\omega}(s))}{\left| \eta_{\omega}'(s) \right|} 
    &+ \beta_2 \int_0^1 \frac{\partial}{\partial n_{\eta_{\omega}(s)}} G_1(\eta_{\omega}(s),\eta_{\omega}(t)) \mu(\eta_{\omega}(t)) dt
    \\ \nonumber
    &+ \alpha_2 \int_0^1 G_1(\eta_{\omega}(s),\eta_{\omega}(t)) \mu(\eta_{\omega}(t)) dt 
	=
	\phi(\eta_{\omega}(s)),&
	\quad
	s \in [0;1],&
\end{align}
which is a Fredholm equation of the second kind.
The pure Neumann problem ($\alpha_2=0$, $\beta_2=1$) yields the similar equation
\begin{align} \label{eq:Neumann_BIE}
	-\frac{1}{2} \frac{\mu(\eta_{\omega}(s))}{\left| \eta_{\omega}'(s) \right|} + \int_0^1 \frac{\partial G_1(\eta_{\omega}(s),\eta_{\omega}(t))}{\partial n_{\eta_{\omega}(s)}} \mu(\eta_{\omega}(t)) dt
	=
	\phi(\eta_{\omega}(s)),
	\quad
	s \in [0;1].
\end{align}

It is well known that equations of the second kind are well-posed.
On the other side, Dirichlet boundary conditions ($\alpha_2=1$, $\beta_2=0$) convert \eqref{eq:BIE} into the Fredholm equation of the first kind
\begin{align} \label{eq:Dirichlet_BIE}
	\int_0^1 G_1(\eta_{\omega}(s),\eta_{\omega}(t)) \mu(\eta_{\omega}(t)) dt
	=
	\phi(\eta_{\omega}(s)),
	\quad
	s \in [0;1],
\end{align}
which is intrinsically ill-posed and must be treated with a special care.
It has been established that equations of the first kind with logarithmically singular kernels admit unique solutions when the conformal radius of the boundary is not equal to one \cite{Yan1988,Cheng1993}.
Therefore, we assume that all boundaries satisfy this condition.

Note that the traditional approach in solving the Dirichlet problems with methods of potential is based on representing the solution in the form of the double-layer potential which results in equations of the second kind.
The use of the single-layer potentials, however, has the advantage of satisfying the governing equation on the deterministic boundary.
As a result, when the length of the boundary $\partial D_1$ is large, the method of Green's potentials can lead to significant computational savings compared to the traditional approaches relying on the discretization of the whole boundary.

Of course, the efficiency of the method of Green's potentials relies on the availability of Green's functions for the specific geometries of the domain. 
Unfortunately, analytical expressions for the Green's functions are known only for very simple domains and the construction of approximations adds an additional level of complexity to the proposed scheme.
However, in the case of uncertain domains, the deterministic complementary problem \eqref{eq:GF_2}-\eqref{eq:GF_BC} has to be solved only once and the value of the Green's function at any field point is then readily available through the simple matrix-vector product which can be done very efficiently.
The implementation aspects of this approach are discussed in succeeding sections.

\section{Discretization scheme}
\label{sec:scheme}

\subsection{Spatial discretization}
\label{sec:spat_discret}

The boundary integral equation \change{\eqref{eq:Neumann_BIE}} has been extensively studied in the literature as the classical equation of potential theory  \cite{Atkinson1991,Atkinson1997,GRAHAM1985,Jeon1997,Sloan2000,Cheng1993,Cheng1995}.
It is an equation of the second kind with a continuous kernel and any approximation technique is expected to work well.
For the sake of completeness, we will use the classical Nystr{\"o}m method. 
We start with the boundary integral operator 
\begin{align*}
	(A\mu)(s) = \int_0^1 K(\eta_{\omega}(s),\eta_{\omega}(t)) \mu(\eta_{\omega}(t)) dt
\end{align*}
with a kernel $K$ and approximate it with the trapezoidal rule on the uniform grid with the step $h~=~1/N$ for some integer $N$
\begin{align} \label{eq:trap_rule}
	(A\mu)(s) \approx (A_h \mu)(s) = h \sum_{k=0}^{N-1} K(\eta_{\omega}(s),\eta_{\omega}(kh)) \mu(\eta_{\omega}(kh)), 
	\qquad 
	s \in [0;1].
\end{align}

After collocating the resulting discrete operator at the quadrature nodes this yields the following numerical scheme for the integral equation in \eqref{eq:Neumann_BIE}
\begin{align}\label{eq:lin_system_neu}
    -\frac{1}{2} \frac{\mu(\eta_{\omega}(k'h))}{\left| \eta_{\omega}'(k'h) \right|} + h \sum_{k=0}^{N-1} \frac{\partial G_1(\eta_{\omega}(k'h),\eta_{\omega}(kh))}{\partial n_{\eta_{\omega}(s)}} \mu(\eta_{\omega}(kh))
    = \phi(\eta_{\omega}(k'h)),& 
    \\ \nonumber
	k' = 0,...,N-1.&
\end{align}
The order of this scheme is determined by the order of the trapezoidal quadrature rule in \eqref{eq:trap_rule}, i.e., it is at least $O(h^2)$ for sufficiently smooth $\mu$ \cite{Atkinson1997}.

The direct extension of the scheme in \eqref{eq:lin_system_neu} to the equations of the first kind is not possible due to the singularity of the kernel. 
In this case, we will use the quadrature technique proposed in \cite{Sloan1992} for the first kind Fredholm equations with logarithmic kernels on closed curves.
It is a fully discrete method of qualocation type based on the composite quadrature rule, i.e., both the integral operator and the Galerkin projection are approximated with suitable quadratures.
In particular, we use the trapezoidal approximation of the integral operator in \eqref{eq:trap_rule} and then project it on the test space $S_h$ of 1-periodic smoothest splines of order $r$ with the discrete inner product
\begin{align*}
	(v,w)_h = Q_h(v {w}),
\end{align*}
where
\begin{align*}
	Q_h g =  h \sum_{k=0}^{N-1} \sum_{j=1}^{J} w_j g((k+\zeta_j)h),
	\qquad
	0 < \zeta_1 < \zeta_2 < ... < \zeta_J < 1
\end{align*}
and
\begin{align*}
	\sum_{j=1}^J w_j = 1, 
	\quad
	w_j > 0,
	\quad
	\text{for } 1 \leq j \leq J.
\end{align*}

The problem now can be formulated as follows: find $\mu_h$ such that 
\begin{align*}
	(A_h \mu_h, \chi)_h = (\phi,\chi)_h,
	\qquad
	\forall \chi \in S_h.
\end{align*}

For $r=2$, the basis $(v_0,...,v_{N-1})$ of $S_h$ is represented by the classical hat functions
\begin{align} \label{eq:basis_fun}
	v_k(s) = 
	\begin{cases}
		1 - |s-kh| / h, & \text{if } |x-kh| \leq h,
		\\
		0, & \text{otherwise}.
	\end{cases}
\end{align}
Given the basis, one can write the discrete formulation of the problem: find $\mu_h$ such that
\begin{align} \label{eq:lin_system}
	\sum_{k=0}^{N-1} a_{l,k} \mu_h(\eta_{\omega}(kh)) = (\phi,v_l)_h,
	\qquad
	l = 0,...,N-1,
\end{align}
where
\begin{align*}
	a_{l,k} =
	h^2 \sum_{k'=0}^{N-1} \sum_{j=1}^{J} w_j G_1\Big( \eta_{\omega}\big((k'+\zeta_j)h\big), \eta_{\omega}(kh) \Big) {v_l}\big((k'+\zeta_j)h\big).
\end{align*}

\begin{figure}[t!]
	\centering
    \includegraphics[width=0.35\textwidth]{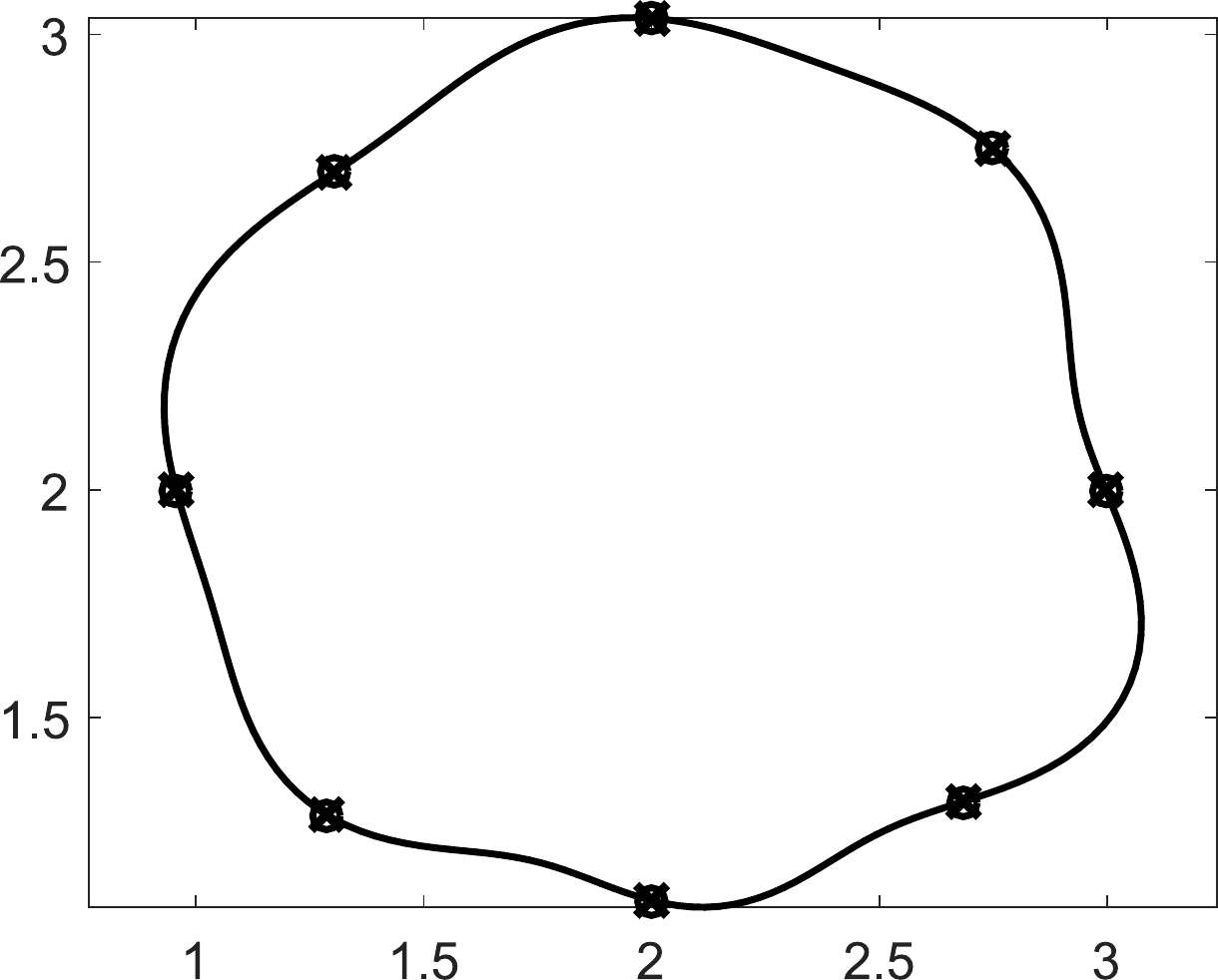}
    \qquad
    \includegraphics[width=0.35\textwidth]{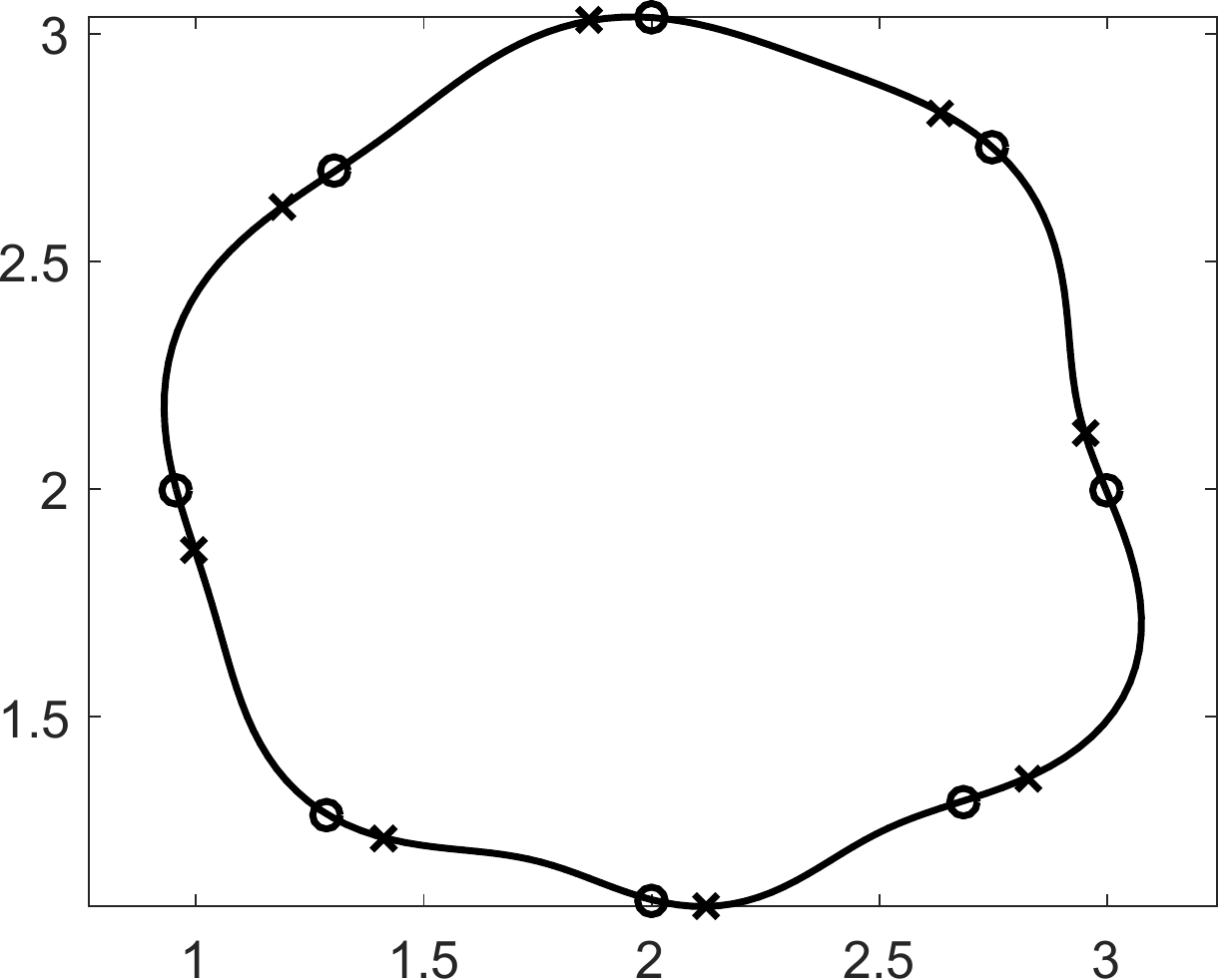}
	\caption{ Trapezoidal rule quadrature points (circles) and collocation points (crosses) for the case of Neumann (left) and Dirichlet (right) boundary conditions. }
	\label{fig:colloc_points}
\end{figure}

It was shown in \cite{Sloan1992} that the following unique choices of the quadrature \change{rules} 
are optimal in terms of stability of the approximation
\begin{equation*}
\begin{minipage}{.45\linewidth}
	\begin{align}\label{eq:colloc_points}
		\nonumber
		J &= 1,  
		\\
		\zeta_1 &= \frac{1}{6},
		\\ \nonumber
		w_1     &= 1.
	\end{align}
\end{minipage}
\begin{minipage}{.1\linewidth}
\qquad
\end{minipage}
\begin{minipage}{.45\linewidth}
    \begin{alignat*}{2}
		\nonumber
		J &= 2, & 
		\\ \tag{\theequation'}
		\zeta_1 &= \frac{1}{6}, \qquad \zeta_2 &= \frac{5}{6},
		\\ \nonumber
		w_1     &= \frac{1}{2}, \qquad     w_2 &= \frac{1}{2},
	\end{alignat*}
\end{minipage}
\end{equation*}
The choice of nodes in (\ref{eq:colloc_points}') gives the $O(h^3)$ order of uniform convergence.
The scheme with nodes in \eqref{eq:colloc_points} has $O(h^2)$ accuracy but the linear system in \eqref{eq:lin_system} has the simpler form
\begin{align} \label{eq:lin_system_2}
	h\sum_{k=0}^{N-1} G_1\Big( \eta_{\omega}\big((k'+\zeta_1)h\big), \eta_{\omega}(kh) \Big) \mu_h(\eta_{\omega}(kh)) = \phi \big( \eta_{\omega}\big((k' + \zeta_1)h \big)\big),&
	\\ \nonumber
	k' = 0,...,N-1.&
\end{align}

The discrete systems in \eqref{eq:lin_system_neu} and \eqref{eq:lin_system_2} have similar form and we use them to solve equations of the second and first kind respectively.
The combined linear system for the problems with both Dirichlet and Neumann boundary conditions can be written as
\begin{align}\label{eq:lin_system_final}
    \mathbf{G} \boldsymbol{\mu} = \boldsymbol{\phi}
\end{align}
where the elements $\mathbf{G}_{k',k}$ of the matrix $\mathbf{G}$ are given by
\begin{align*}
    \mathbf{G}_{k',k} = 
    \left\{
		\begin{array}{l l}
			h G_1\Big( \eta_{\omega}\big((k'+\zeta_1)h\big), \eta_{\omega}(kh) \Big) &\qquad \text{if  } \eta_{\omega}(k'h)\in\mathcal{D},
			\\[1em]
			-\frac{1}{2} \frac{\mu(\eta_{\omega}(k'h))}{\left| \eta_{\omega}'(k'h) \right|} \mathcal{I}_{k'-k} + h \frac{\partial G_1(\eta_{\omega}(k'h),\eta_{\omega}(kh))}{\partial n_{\eta_{\omega}(s)}} &\qquad \text{if  } \eta_{\omega}(k'h)\in\mathcal{N},
		\end{array}
	\right.
\end{align*}
$\mathcal{I}_{k'-k}$ is the indicator function and $\mathcal{D}$, $\mathcal{N}$ denote the collections of nodes with Dirichlet and Neumann boundary conditions respectively.
The quadrature and collocation points for the case of the Neumann and Dirichlet boundary conditions are depicted in Figure \ref{fig:colloc_points}.

After the density of the potential is determined from \change{\eqref{eq:lin_system_final}}, one can calculate the solution $u_2(x,\omega)$ at any field point by evaluating the integral with the trapezoidal quadrature rule \eqref{eq:trap_rule} as
\begin{align*}
	u_2(x,\omega) 
	\approx h \sum_{k=0}^{N-1} G_1(x,\eta_{\omega}(kh)) \mu(\eta_{\omega}(kh)).
\end{align*}

\subsection{Evaluation of Green's functions for arbitrary domains}
\label{sec:num_GF}

Recall the definition of the Green's function for the Laplace operator as a solution of the following boundary value problem
\begin{alignat*}{3}
    \tag{\ref{eq:GF_2}}
	-\nabla^2 G_1(x,\xi)                                                     & = \delta(\xi)  &\qquad &\text{for } x, \xi \in D_1,
	\\\tag{\ref{eq:GF_BC}}
	\alpha_1 G_1(x,\xi) + \beta_1 \frac{\partial G_1(x,\xi)}{\partial n_x}   & = 0  &\qquad &\text{for } x \in \partial D_1.
\end{alignat*}
It was mentioned previously that the proposed numerical technique relies heavily on the ability to evaluate Green's functions for the domains of arbitrary shapes.
We outline here several methods which allow to solve \eqref{eq:GF_2}-\eqref{eq:GF_BC} in a computationally attractive way.

\subsubsection{Analytical Green's functions.}

In certain cases, it is possible to solve \eqref{eq:GF_2}-\eqref{eq:GF_BC} analytically.
The \textbf{\textit{fundamental solution}} of the 2-D Laplace equation is one of such examples of an exceptional importance. 
It satisfies the equation \eqref{eq:GF_2} in the entire space and has the form
\begin{align} \label{eq:fund_sol}
	\mathfrak{G}(x,\xi) = -\frac{1}{2 \pi} \ln r,
	\qquad
	r = \sqrt{ (x_1-\xi_1)^2 + (x_2-\xi_2)^2 }.
\end{align}

The Green's function to \eqref{eq:GF_2}-\eqref{eq:GF_BC} can then be written as a sum of \eqref{eq:fund_sol} and a ``corrector" function aiming to satisfy the boundary condition in \eqref{eq:GF_BC}
\begin{align} \label{eq:GF_general}
	G_1(x,\xi) = \mathfrak{G}(x,\xi) + \psi(x, \xi).
\end{align}
Defined in this way, the corrector function is called the regular component of the Green's function.
It solves the following complementary problem
\begin{alignat}{3} \label{eq:corrector}
	-\nabla^2 \psi(x,\xi)                                                   & = 0  &\qquad &\text{for } x, \xi \in D_1,
	\\ \label{eq:corrector_BC}
	\alpha_1 \psi(x,\xi) + \beta_1 \frac{\partial \psi(x,\xi)}{\partial n_x}   & = -\alpha_1 \mathfrak{G}(x,\xi) - \beta_1 \frac{\partial \mathfrak{G}(x,\xi)}{\partial n_x}  &\qquad &\text{for } x \in \partial D_1.
\end{alignat}
For simple geometries possessing certain symmetry properties, the problem in \eqref{eq:corrector}-\eqref{eq:corrector_BC} admits the closed form solution which can be constructed via the method of images.
Two classical examples of such Green's functions are given below.

\textit{\textbf{Dirichlet problem in the upper half-plane}} $D_1(x_1,x_2)=\{x_1,x_2\geq0\}$
\begin{align*}
	G_1(x,\xi) = \frac{1}{2 \pi} \ln \sqrt{ \frac{ (x_1-\xi_1)^2 + (x_2+\xi_2)^2 }{ (x_1-\xi_1)^2 + (x_2-\xi_2)^2 } }.
\end{align*}

\textit{\textbf{Dirichlet problem in the disk}} $D_1(r,\varphi)=\{0 \leq r < a, 0 \leq \varphi \leq 2 \pi \}$
\begin{align*}
	&G_1 \big( x, \xi \big) = \frac{1}{4 \pi} \ln \frac{a^4 - 2a^2r\rho  \cos(\varphi- \varsigma)+r^2 \rho^2}{a^2(r^2-2r\rho \cos(\varphi- \varsigma)+\rho^2)},
	\quad
	\\[1em]
	&(x_1,x_2) = r ( \cos(\varphi), \sin(\varphi) ),
	\quad
	(\xi_1,\xi_2) = \rho ( \cos(\varsigma), \sin(\varsigma) ).
\end{align*}

More examples can be found, for instance, in \cite{Duffy2001,Melnikov2011}.
Additionally, the infinite product representation of Green's functions arising from applying the method of images was discussed in \cite{Melnikov2012,Reshniak2013}.

\subsubsection{Direct numerical approximation of Green's functions.}

By its definition, Green's function of the boundary value problem is the inverse of the corresponding differential operator.
Numerical methods can be viewed as implicitly approximating such inverse operators.
For instance, it was shown in \cite{Tottenham1970} that the finite element solution $u_h(x)$ of the boundary value problem with homogeneous boundary conditions has the form
\begin{align*}
	u_h(x) = \int_D G_h(x,\xi) f(\xi) d \xi,
\end{align*}
where $G_h(x,\xi)$ is the FE-Green's function, i.e., the projection of the exact Green's function on the finite element space $V_h$.
Solving for $G_h(x,\xi)$ yields

\begin{Theorem}[\cite{Hartmann2013}] 
Let $\mathbf K$ be the stiffness matrix of the linear system arising after the finite element disretization.
The FE-Green's function has the form
\begin{align} \label{eq:FE_GF}
	G_h(x,\xi) = \mathbf{v}(x)^T \mathbf K^{-1} \mathbf{v}(\xi),
\end{align}
where $\mathbf{v}(x)=(v_1(x),...,v_M(x))$ are the basis functions of the FE-space $V_h$.
\end{Theorem}

For the spectral finite element method, the formula \eqref{eq:FE_GF} simplifies to
\begin{align} \label{eq:spectral_GF}
	G_h(x,\xi) = \mathbf{v}(x)^T \mathbf\Lambda^{-1} \mathbf{v}(\xi) = \sum_{i=1}^M \frac{ v_i(x) v_i(\xi) }{\lambda_i},
\end{align}
where $(v_1,...,v_M) \in V_h$ are the $M$ leading eigenfunctions of the differential operator and $\mathbf\Lambda$ is the diagonal matrix with the corresponding eigenvalues.
One can immediately recognize in \eqref{eq:spectral_GF} the truncation of the classical eigenfunction representation of the Green's function which is guaranteed to exist by the Mercer's theorem \cite{Fasshauer2016}.
Spectral representations of the Green's functions can be found, for instance, in \cite{Duffy2001}.
As an example, we provide the Green's function for the 

\textit{\textbf{Dirichlet problem in the rectangle}} $D_1(x_1,x_2)=\{0 \leq x_1 \leq a, 0 \leq x_2 \leq b \}$
\begin{align} \label{eq:spectral_GF_rect}
	&G_1 \big( x, \xi \big) = 
	4ab \sum_{n=1}^{\infty} \sum_{m=1}^{\infty} \frac{ \displaystyle{ \sin \left( \frac{m\pi x_1 }{a} \right) \sin \left( \frac{n\pi x_2 }{b} \right) \sin \left( \frac{m\pi \xi_1 }{a} \right) \sin \left( \frac{n\pi \xi_2 }{b} \right)  } }{ n^2 \pi^2 a^2 + m^2 \pi^2 b^2 }.
\end{align}

\subsubsection{Numerical approximation of the regular part of Green's functions.}

Green's functions generally do not belong to the function spaces approximated by the $span(v_1(x),...,v_M(x))$ resulting in the very slow convergence of the representations in \eqref{eq:FE_GF}-\eqref{eq:spectral_GF}.
For instance, solutions to the Poisson equation are usually constructed in $H^1(D_1)$ but the solution to the problem in \eqref{eq:GF_2}-\eqref{eq:GF_BC} is not in $H^1(D_1)$ since the delta function $\delta(\xi) \notin H^{-1}$ for $d\geq2$, where $d$ is the physical dimension of the problem.
Analogously, analytical expansions like \eqref{eq:spectral_GF_rect} do not have uniform error estimates which seriously limits their immediate practical utilization.

In certain cases, one can obtain uniformly convergent spectral representations by partial summation of the series leading to the explicit extraction of the singularity.
Justification of this approach with practical examples can be found in \cite{Melnikov1998}.
For instance, the series in \eqref{eq:spectral_GF_rect} is transformed to the following form
\begin{align} \label{eq:spectral_GF_rect_summed}
	G_1 \big( x, \xi \big) 
	&= 
	\frac{1}{2\pi} \ln \left[ \frac{E(z-\zeta^{*})E(z+\zeta^{*})E(z_1+\zeta_1^{*})E(z_2+\zeta_2^{*})}{E(z-\zeta)E(z+\zeta)E(z_1+\zeta_1)E(z_2+\zeta_2)} \right]
	\\ \nonumber
	&-
	\frac{2}{b} \sum_{n=1}^{\infty} S_n(x_1,\xi_1) \sin\left(\nu\xi_2\right) \sin\left(\nu x_2\right),
\end{align}
where $\nu = n\pi / b$, $z=x_1+ix_2$, $z_1=(x_1+a)+i x_2$, $z_2=(x_1-a)+i x_2$, $\zeta=\xi_1+i \xi_2$, $\zeta_1=(\xi_1+a)+i \xi_2$, $\zeta_2=(\xi_1-a)+i \xi_2$, $\zeta_1^{*}=(\xi_1+a)-i \xi_2$, $\zeta_2^{*}=(\xi_1-a)-i \xi_2$, $E(z)=\left|e^{\pi z/b}-1\right|$ and
$$
S_n(x_1,\xi_1) = \frac{e^{\nu x_1}\sinh(\nu(\xi_1-a))-e^{-\nu x_1}\sinh(\nu(\xi_1+a))}{2\nu e^{2\nu a}\sinh(\nu a)}.
$$
The remainder term $R_M(x,\xi)$ of the $M$-term truncation of the expansion in \eqref{eq:spectral_GF_rect_summed} has the upper bound 
$$
|R_M(x,\xi)| \leq \frac{b}{2\pi} \left( \ln\left(1-e^{-\pi a/b}\right) - \sum_{n=1}^N \frac{e^{-n\pi a/b}}{n}  \right)
$$
which reveals the extremely high rate of convergence.

Similarly, the corrector function $\psi(x,\xi)$  in \eqref{eq:GF_general} is harmonic everywhere in $D_1$ and thus can be efficiently approximated with any conventional numerical method.
For instance, the finite element approximation has the form
\begin{align}\label{eq:GF_FE_reg_part}
	\psi_h(x,\xi) = \mathbf{v}(x)^T \mathbf K^{-1} \boldsymbol{\mathfrak{g}}(\xi),
\end{align}
where $\mathbf K$ is the same stiffness matrix as in \eqref{eq:FE_GF} and $\boldsymbol{\mathfrak{g}}(\xi)$ encodes the \change{negative} trace of the fundamental solution on the boundary.

One can also construct the regular part of the Green's function in the form of the single layer potential
\begin{align}\label{eq:GF_reg_part}
	\psi(x,\change{\xi}) 
	= \int_{\partial D_1} \mathfrak{G}(x,y) \nu_1(y,\change{\xi}) d l(y) 
	= \int_0^1 \mathfrak{G}(x,\eta_1(t)) \mu^{\psi}(\eta_1(t),\change{\xi}) dt,
\end{align}
where $\mathfrak{G}(x,y)$ is the fundamental solution of the differential operator, $\mu^{\psi}(\eta_1(t),\change{\xi}) = \nu_1(\eta_1(t),\change{\xi}) \left| \eta_1'(t) \right|$ and $\eta_1(t)$ defines a parameterization of the boundary $\partial D_1$. 
It is natural to build the approximate solution of the above equation with the same method used for the approximation of the original integral equation, e.g., with the scheme given in section \ref{sec:spat_discret}.
For this, note that similarly to \eqref{eq:GF_FE_reg_part}, the density $\mu^{\psi}(\eta_1(t),\change{\xi})$ encodes the negative trace of the fundamental solution $\mathfrak{G}(\cdot,\xi)$ on the boundary $\partial D_1$.

\begin{figure}[!t]
    \centering
    \includegraphics[width=0.49\textwidth]{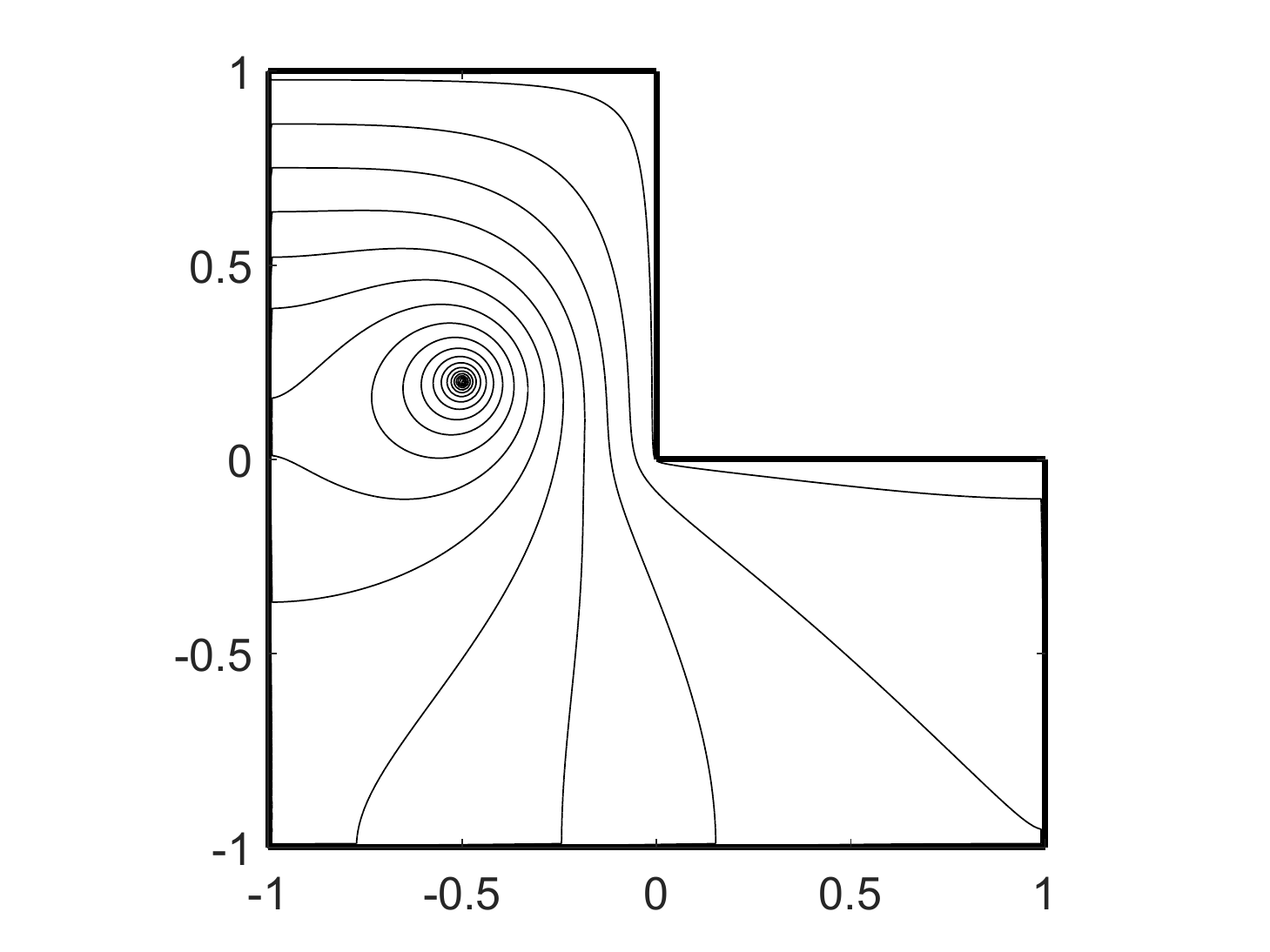}
    \includegraphics[width=0.49\textwidth]{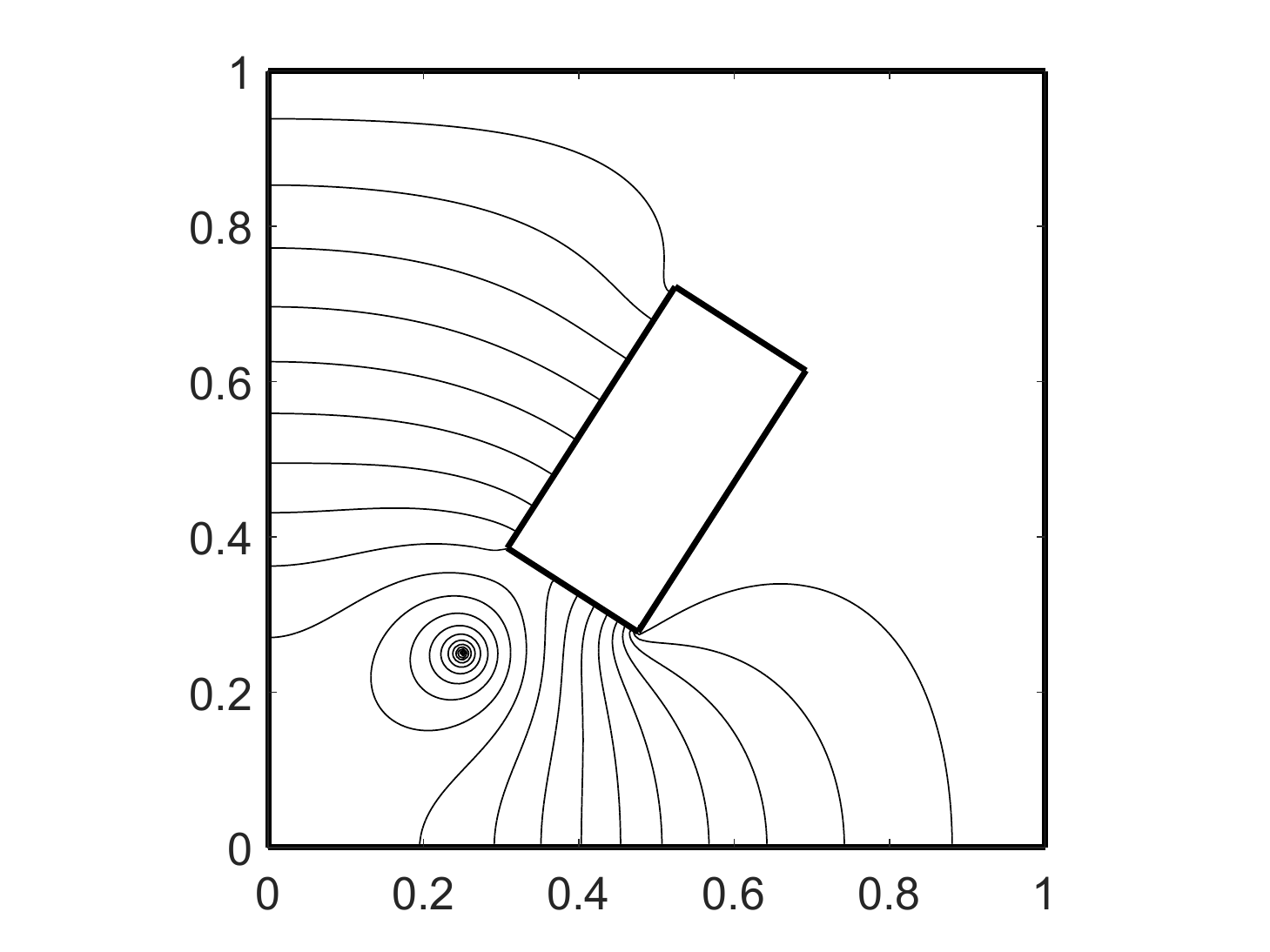}
	\caption{ Numerical Green's functions in the L-shaped and multiconnected domains with a combination of Dirichlet and Neumann boundary conditions.}
	\label{fig:GF_num}
\end{figure}

Denote by $\boldsymbol{\mathfrak{G}}$ the marix of the linear system in \eqref{eq:lin_system_final} for the problem formulated over the whole boundary $D=D_1\cup D_2$ and using fundamental solution $\mathfrak{G}(x,\xi)$ instead of the Green's function $G_1(x,\xi)$.
By taking $\boldsymbol\mu^{\psi}(\xi) = \big( \mu^{\psi}_{0}(\xi), ..., \mu^{\psi}_{N-1}(\xi) \big) $ with $\mu^{\psi}_{k}(\xi) = \mu^{\psi}(\eta_1(kh),\xi)$, we get
\begin{align*}
    \boldsymbol\mu^{\psi}(\xi) = \boldsymbol{\mathfrak{G}}_{11}^{-1} \mathfrak{g}(\xi),
\end{align*}
where $\boldsymbol{\mathfrak{G}}_{11}$ is a $N_1\times N_1$ submatrix of $\boldsymbol{\mathfrak{G}}$ with both $\eta_{\omega}(k'h)\in D_1$ and  $\eta_{\omega}(kh)\in D_1$.
The matrix $\boldsymbol{\mathfrak{G}}_{11}$ needs to be inverted only once and the approximation of the regular part of the Green's function at any point is then readily available as
\begin{align*}
    &\psi_h(x,\xi) 
	= h \sum_{k=0}^{N-1} \mathfrak{G}(x,\eta_1(kh)) \mu^{\psi}(\eta_1(kh),\xi).
\end{align*}
It is easy to verify that in this case the matrix $\mathbf{G}\in\mathbb{R}^{N_2\times N_2}$ of the original linear system in \eqref{eq:lin_system_final} takes the form of the Schur complement of $\boldsymbol{\mathfrak{G}}_{11}$
\begin{align}\label{eq:GF_num_reg_part}
    \boldsymbol{G}
	= \boldsymbol{\mathfrak{G}}_{22} - \boldsymbol{\mathfrak{G}}_{21} \boldsymbol{\mathfrak{G}}_{11}^{-1} \boldsymbol{\mathfrak{G}}_{12},
\end{align}
where by analogy with $\boldsymbol{\mathfrak{G}}_{11}$, $\boldsymbol{\mathfrak{G}}_{22}$ is a $N_2\times N_2$ submatrix of $\boldsymbol{\mathfrak{G}}$ with $\eta_{\omega}(k'h)\in D_2$ and $\eta_{\omega}(kh)\in D_2$, and similarly for $\boldsymbol{\mathfrak{G}}_{21}\in\mathbb{R}^{N_2\times N_1}$ and $\boldsymbol{\mathfrak{G}}_{12}\in\mathbb{R}^{N_1\times N_2}$.
We will use this approach in the subsequent sections.
Two examples of the approximate Green's function obtained in this way are given in Figure \ref{fig:GF_num}.

\subsection{Statistical discretization}

Let $u(x,\omega)$ be the solution of the random partial differential equation and denote by $f(\omega):=f(u(x,\omega))$ some functional of $u(x,\omega)$.
The expected value of $f(\omega)$ can be approximated by the Monte-Carlo (MC) estimator of the form
\begin{align*}
	\mathbb{E} [ f(\omega) ] 
	\approx 
	\EMC{f} = \frac{1}{M} \sum_{m=1}^M f^m({x}),
\end{align*}
where the deterministic functions $f^m$ represent i.i.d. realizations of $f(\omega)$. 

Monte Carlo method is a purely statistical technique which ignores any information about regularity of the functions in the physical space.
The multilevel Monte Carlo method can exploit some of this information by relating the sampling error of the estimator to the convergence properties of the spatial discretization \cite{Giles2008}. 
For instance, consider a hierarchical family of nested discretizations with step sizes
\begin{align*}
	h_L < h_{L-1} < ... < h_l < ... < h_0, \qquad h_l = q^{-l} h_0,
\end{align*}
where $q \in \mathbb{N} \setminus 1$ is a refinement parameter.
Denote by $u_l(x,\omega)$ the approximation of $u(x,\omega)$ at the level $l$ and let $f_{l} (\omega) = f(u_l ({x},\omega))$. Then the solution at the finest discretization level $L$ is given as the telescoping sum
\begin{align*}
	f_{L} (\omega) = f_{0} (\omega) + \sum_{l=1}^L \Big( f_{l} (\omega) - f_{l-1} (\omega) \Big).
\end{align*}
Taking advantage of the linearity of the expectation, we obtain
\begin{align*}
	\mathbb{E} \Big[ f_{L} (\omega) \Big] = \mathbb{E} \Big[ f_{0} (\omega) \Big] + \sum_{l=1}^L \mathbb{E} \Big[ f_{l} (\omega) - f_{l-1} (\omega) \Big].
\end{align*}

By setting
$	
	\Delta_{l} ( \omega )
	= f_{l} ( \omega ) - f_{l-1} ( \omega )
	, \; l=0,...,L,
$
the above expression yields the multilevel Monte-Carlo estimator
\begin{align}\label{eq:MLMC_estimator}
	\E{ f({x},\omega) }
	\approx
	\EML{f_L}
	= \sum_{l=0}^L \EMC{\Delta_{l}}
	= \sum_{l=0}^L M_l^{-1} \sum_{m_l=1}^{M_l} \Delta_{l}^{m_l} 
	,
\end{align}
where $\Delta_{0} ( \omega ) = f_{0} ( \omega )$ and $M_l$ is the number of random samples at the level $l$.

\section{Complexity analysis}
\label{sec:complexity}

\subsection{Asymptotical complexity of the MLMC estimator}

Consider the cost function of the MLMC estimator
\begin{align*}
	C_{\mathbb{ML}} 
	&= \sum_{l=0}^L M_l {C}_l 
	,
\end{align*}
where $M_l$ is the number of samples at the level $l$  and $C_l$ is the cost of generating the single realization of $f_l(\omega)$.
The following theorem shows that the above cost can be optimized by appropriately balancing the errors of the Monte Carlo estimators at different levels.

\begin{Theorem}[\cite{Giles2008}]\label{the:I_MLMC_cost}
If there exist independent estimators $\EMC{\Delta_{l}}$ based on $M_l$ Monte Carlo samples, each with expected cost $C_l$ and variance $V_l$, and positive constants $\alpha$, $\beta$, $\rho$ such that $\min(\beta, \rho) \leq 2\alpha$ and
\begin{enumerate}
	\item $\Big| \E{f_l(\omega) - f \big( u(x,\omega \big)} \Big| \lesssim h_l^{\alpha}$,
	\item $V_l \lesssim h_l^{\beta}$,
	\item $C_l \lesssim h_l^{-\rho}$,
\end{enumerate}
then for any $\epsilon < e^{-1}$ there are values $L$ and $M_l$ for which the multilevel estimator \eqref{eq:MLMC_estimator} has a mean-square-error with bound 
\begin{align*}
	\E{ \Big( \EML{f_L(\omega)} - \E{f\big(u(x,\omega)\big)} \Big)^2 }
	< \epsilon^2
\end{align*}
with a computational complexity $C_{\mathbb{ML}}$ with bound
	\begin{align}\label{eq:I_MLMC_cost}
		C_{\mathbb{ML}}
		\lesssim 
		\begin{cases}
			\epsilon^{-2}  & \text{if } \rho - \beta < 0, \\
			\epsilon^{-2} \big| \ln \epsilon \big|^2  & \text{if } \rho - \beta = 0, \\
			\epsilon^{-2-\frac{\rho-\beta}{\alpha}} & \text{if } \rho - \beta > 0.
		\end{cases}
	\end{align}
\end{Theorem}

\textbf{Corollary.}
According to Theorem \ref{the:I_MLMC_cost}, the optimal $\epsilon^{-2}$ complexity of the MLMC estimator is achieved when $\rho<\beta$, i.e., when the decay of the variance $V_l$ of the level corrections is faster than the growth of the computational cost $C_l$.
For Lipschitz continuous $f$ and given the convergence rate $\alpha$ of the spatial discretization scheme, the variance $V_l$ decays with order $\beta=2\alpha$.
In practice, $\rho\in(2,3]$ for different methods yielding $\alpha>\frac{\rho}{2}\in (1,1.5]$ as an optimality condition which is easily satisfied with most practical schemes including the one utilized in this paper.
However, the complexity estimate in \eqref{eq:I_MLMC_cost} has a silent constant which determines the actual relative performace of different algorithms with the same asymptotic complexity.
In the following sections, we provide the complexity analysis for the proposed numerical scheme and use it to justify our choice of the method of Green's potentials.

\subsection{Error component analysis}
Consider the following decomposition of the total error of the MLMC estimator
\begin{align} \label{eq:error_decomp}
	&\norm{ \EML{\tilde{u}_L} - \E{ u } }{L^2(\Omega)}
	\leq \underbrace{ \norm{ \E{ \tilde{u}_L - u }     }{H^1(D)}  }_{\text{I := Discretization error}}
	+     \underbrace{ \norm{ \dEML{ \tilde{u}_L} }{L^2(\Omega)} }_{\text{II := Sampling error}}
	, 
\end{align}
where $u$ and $\tilde{u}_L$ are the exact and the approximate values of the potential \eqref{eq:potential}. 
The error components $I$ and $II$ correspond to the spatial approximation error and the sampling error.
To achieve the desired accuracy $\epsilon$, it is sufficient to balance the total error between these two components in the following way
\begin{align} \label{eq:MLMC_tol}
	\norm{ \EML{\tilde{u}_L} - \E{ u } }{L^2(\Omega)}
	&\leq \epsilon_I + \epsilon_{II} 
	= \epsilon
	.
\end{align}

The error analysis for each of the components is provided below.

\subsubsection{Spatial discretization error.}

Jensen's inequality gives
\begin{align}\label{eq:eps_I}
	\norm{ \E{ u_L({x},\omega) - u ({x},\omega)} }{H^1(D)}
	&\leq
	\E{ \norm{ u_L({x},\omega) - u ({x},\omega) }{H^1(D)} }
	= \epsilon_{I}
\end{align}
and the estimate of the error component $I$ can be derived from the convergence properties of the spatial discretization scheme.
The order of such scheme is in turn controlled by the convergence properties of the boundary integral equation solver and the numerical error of the single-layer potential evaluation.

For instance, it was shown in \cite{Sloan1992} that the qualocation-type discretization \eqref{eq:lin_system_2} of the first-kind integral equation in \eqref{eq:Dirichlet_BIE} admits the following estimate for the approximation error of $\mu \in H^t([0,1])$
\begin{align}\label{eq:BIE_error}
\norm{\mu_{l}-\mu}{H^s([0,1])} \leq c h_l^{t-s} \norm{\mu}{H^t([0,1])}
\end{align}
provided that $s>\frac{1}{2}$, $s+\frac{1}{2}< t \leq s+\alpha$ and the right hand side of the integral equation is continuous and 1-periodic.
In the case of optimal regularity, i.e., for $\mu \in H^t([0,1])$ with $t>\alpha+1/2$, the following error bound is valid
\begin{align}\label{eq:spat_order}
\sup_{t \in [0,1]}| \mu_l(t)-\mu(t) | = \norm{\mu_l-\mu}{L^{\infty}([0,1])} \leq c h_l^{\alpha} \norm{\mu}{H^t([0,1])}
\end{align}
due to the embedding of $H^s$ ($s>1/2$) in $C_p$, the space of 1-periodic continuous functions.
The order of convergence is $\alpha=2$ for the scheme in \eqref{eq:lin_system_2} and hence the method is $O(h_l^2)$ accurate.

The error of the Nystr\"om method in \eqref{eq:lin_system_neu} for the integral equations of the second kind is given by the error of the corresponding numerical quadrature rule.
Since the trapezoidal rule is exteremely well-behaved for smooth periodic functions, the order of the Nystr\"om-trapezoidal method is also at least $O(h_l^2)$ \cite{Atkinson1997}.

To determine the overall order of the spatial discretization at the level $l$, consider the error in approximating the single-layer potential \eqref{eq:potential}
 \begin{align*}
|\tilde{u}_l(x) - u(x)|
&\leq |u_l(x) - u(x)| + |\tilde{u}_l(x) - u_l(x)|,
\end{align*}
where the first term on the right side gives the error of the numerical scheme with the exact Green's function and the second term is the error due to the approximation of the Green's kernel itself.
We get for the first component that
\begin{align*}
u_l(x) - u(x)
&= h_l \sum_{k=0}^{N_l-1} G_1(x,\eta_{\omega}(kh)) \mu(kh_l) - \int_0^1 G_1(x,\eta_{\omega}(t)) \mu(t) dt
\\ \nonumber
&=  \int_0^1 G_1(x,\eta_{\omega}(t)) ( \mu_l(t) - \mu(t) ) dt  +  R(x),
\end{align*}
where $R(x)$ is the error of the trapezoidal rule in \eqref{eq:trap_rule}.
For a periodic function $\mu\in H^t$, $R(x)$ is also $O(h^t)$ \cite[Lemma~A1]{Sloan1992}, and since $t>\alpha$ in \eqref{eq:spat_order}, one gets
\begin{align}\label{eq:pot_error}
| u_l(x) - u(x) |
&\leq \left| \int_0^1 G_1(x,\eta_{\omega}(t)) ( \mu_l(t) - \mu(t) ) dt \right| + O(h_l^t)
\\ \nonumber
&\leq
\norm{\mu_l-\mu}{L^{\infty}([0,1])} \int_0^1 | G_1(x,\eta_{\omega}(t)) |  dt + O(h_l^t)
= c \norm{\mu_l-\mu}{L^{\infty}([0,1])}.
\end{align}
Similarly, the Green's kernel is a harmonic and thus analytic function at any internal point of the domain which gives the estimate for the error in the approximation of the derivatives of the potential
\begin{align}\label{eq:pot_deriv_error}
| u_l^{(i)}(x) - u^{(i)}(x) |
&\leq \left| \int_0^1 G_1^{(i)}(x,\eta_{\omega}(t)) \big( \mu_l(t) - \mu(t) \big) dt \right| + O(h_l^t)
\\ \nonumber
&\leq \norm{\mu_l-\mu}{L^{\infty}([0,1])} \int_0^1 |  G_1^{(i)}(x,\eta_{\omega}(t)) |  dt + O(h_l^t)
= c \norm{\mu_l-\mu}{L^{\infty}([0,1])},
\end{align}
where $i=(i_1, i_2)$ is a multi-index and $f^{(i)} = \displaystyle{ \frac{\partial^{|i|} f}{\partial x_1^{i_1} \partial x_2^{i_2}} }$.

Now consider the approximation of the Green's kernel.
From \eqref{eq:GF_general} and \eqref{eq:GF_reg_part}, we have 
\begin{align*}
	G_1(x,\xi) = \mathfrak{G}(x,\xi) +  \int_0^1 \mathfrak{G}(x,\eta_{1}(t)) \mu^{\psi}(\eta_{1}(t),\xi) dt.
\end{align*}
By analogy with \eqref{eq:pot_error}, there holds the error estimate
\begin{align}\label{eq:GF_error}
\left| \tilde{G}_1(x,\cdot) - G_1(x,\cdot) \right|
= \left| \psi_{l}(x,\cdot) - \psi(x,\cdot) \right|
= c \norm{\mu_{l}^{\psi}-\mu^{\psi}}{L^{\infty}([0,1])}.
\end{align}

By combining \eqref{eq:spat_order}, \eqref{eq:pot_error} and \eqref{eq:GF_error}, one gets the error bound
\begin{align*}
| \tilde{u}_l(x) - u_l(x) |
&\leq h_l \cdot \sup_{k} |\mu(kh_l)| \cdot \sum_{k=0}^{N_l-1} \left| \tilde{G}_1(x,\xi(kh_l)) - G_1(x,\xi(kh_l)) \right|
= c h_l^{\alpha}.
\end{align*}

The estimate for the first error component follows trivially from \eqref{eq:pot_deriv_error} as
\begin{align}\label{eq:H_s_pot_err}
I &:= \E{\norm{ u_l(x,\omega) - u(x,\omega)}{H^1(D)}}
\leq c h_l^{\alpha} = \epsilon_I.
\end{align}
Hence, the condition in \eqref{eq:eps_I} on the spatial discretization error can be satisfied by choosing the number of levels according to
\begin{align*}
	h_L = q^{-L} h_0
	\quad
	\to
	\quad
	L = \left\lceil \log_q \left( h_0 (c_1 \epsilon_I^{-1})^{1/\alpha} \right) \right\rceil
	\leq c + \log_q \left( h_0 \epsilon_I^{-1/\alpha} \right).
\end{align*}

\textbf{Remark.} 
We used the embedding argument in \eqref{eq:spat_order}, the analyticity of the Green's kernel away from the boundary and the high accuracy of the trapezoidal rule for smooth periodic functions to find the errors \eqref{eq:pot_error}-\eqref{eq:pot_deriv_error} of the PDE solution $u(x)$ and its derivatives $u^{(i)}(x)$.
For Galerkin methods, similar estimates can be determined directly from the energy errors of the form \eqref{eq:BIE_error} using the Aubin-Nitsche-duality argument \cite{sauter2010boundary,steinbach2007numerical}.
Any such method to solve boundary integral equations can be trivially applied instead of the one utilized in this paper which might be of interest for the problems in higher dimensions or in regions with piecewise smooth boundaries.
In this case, the results of the MLMC theory remain valid but with possibly different values for the rates $\alpha$, $\beta$, $\rho$ in conditions $(1)-(3)$ of Theorem \ref{the:I_MLMC_cost}. 
For planar regions with analytic contours, the high accuracy of the quadrature based methods combined with the trivial construction of their matrices motivated our choice of the solvers in \eqref{eq:lin_system_neu}-\eqref{eq:lin_system_2}.

\subsubsection{Sampling error.}
	From the definition of the norm in \eqref{eq:rand_norm}, we obtain
	\begin{align*}
		&\norm{ \dEMC{u} }{L^2(\Omega)}^2
		=
		\E{ \norm{ \dEMC{u} }{H^1(D)}^2 }
		\\
		&\qquad
		=
		\frac{1}{M^2} \sum_{i=0}^s \int_D \sum_{m=1}^M \mathbb{V}ar \Big( u^{(i),m} \Big) d {x}
		+ \frac{1}{M^2} \sum_{i=0}^s \int_D \sum_{m=1}^M \sum_{\substack{m'=1 \\ m \neq m'}}^M \mathbb{C}ov \Big( u^{(i),m},u^{(i),m'} \Big) d {x}.
	\end{align*}
	By virtue of the independence of i.i.d. samples $u^m$, we have $\mathbb{C}ov \Big( u^{(i),m},u^{(i),m'} \Big)~=~0$ and
	\begin{align}\label{eq:MC_error}
		\norm{ \dEMC{u} }{L^2(\Omega)}
		=
		\sqrt{ \frac{\overline{V} ( u )}{M}   },
	\end{align}	
	where $\overline{V} ( u ) = \E{ \norm{ u - \E{u} }{H^1(D)}^2 }$.

Due to the independence of the MC estimators at each level, one gets the error of the MLMC estimator as follows
\begin{align*}
	\norm{ \dEML{ u_{L}} }{L^2(\Omega)}^2
	&= \norm{ \E{ u_{L}} - \sum_{l=0}^L \EMC{ u_{l} - u_{{l-1}} } }{L^2(\Omega)}^2
	\\ \nonumber
	& = \sum_{l=0}^L \E{ \norm{ \EMC{ \Delta_{l} } - \E{ \Delta_{l} } }{H^1(D)}^2 }
	 =
	\sum_{l=0}^L \frac{\overline{V}_l}{M_l},
\end{align*}
where $\overline{V}_l=\overline{V}(\Delta_{l})$ and $\Delta_{l}=u_{l} - u_{{l-1}}$.

Thus, the sampling error can be estimated as
\begin{align} \label{eq:sampling_err}
	II^2 &:=
	\norm{ \dEML{ u_{L}} }{L^2(\Omega)}^2
	=
	\sum_{l=0}^L \frac{\overline{V}_l}{M_l}
	= \epsilon_{II}^2.
\end{align}
The above condition can be satisfied by taking
\begin{align*}
	M_l = \frac{ \overline{V}_l }{ a_l \epsilon_{II}^2 }
	\qquad
	\text{ and }
	\qquad
	\sum_{l=0}^L a_l = 1,
\end{align*}
where the coefficients $a_l$ are the weights assigning certain part of the sampling error to each level.
It was shown in \cite{Giles2008} that the optimal cost in Theorem \ref{the:I_MLMC_cost} is achieved with the following choice
\begin{align}\label{eq:num_samples}
	a_l 
	= 
	\frac{({C_l} \overline{V}_l )^{1/2} }{ \displaystyle{\sum_{k=0}^L ( {C_k} \overline{V}_k )^{1/2} }}
	\quad
	\to
	\quad
	M_l = \epsilon_{II}^{-2} \left( \frac{ \overline{V}_l }{ C_l  } \right)^{1/2} \sum_{k=0}^L ( {C_k} \overline{V}_k )^{1/2}.
\end{align}

Note that, in view of \eqref{eq:H_s_pot_err}, one has $\overline{V}_l = O(h_l^{2 \alpha})$, i.e., $\beta = 2 \alpha \geq 4$ in Theorem~\ref{the:I_MLMC_cost}.

\subsubsection{Complexity of the method of Green's potentials with different kernels}

In this section, we study and compare the complexity of the method of Green's potentials for three different choices of the kernel function: fundamental solution, analytical Green's kernel and approximate Green's kernel.

\subsubsection{Fundamental solution}
\label{sec:Schur_complexity}

Here we provide the complexity analysis of the scheme in section \ref{sec:spat_discret} with the boundary integral equations formulated on the whole boundary $\partial D = \partial D_1~\cup~\partial D_2$.
Consider the overall cost of the MLMC estimator with $M_l$ as in \eqref{eq:num_samples}
\begin{align*}
	C_{\mathbb{ML}} 
	&= \sum_{l=0}^L \Big \lceil M_l \Big \rceil {C}_l 
	\leq 
	\sum_{l=0}^L {C}_l + \epsilon_{II}^{-2} \left( \sum_{l=0}^L ({C_l} \overline{V}_l )^{1/2} \right)^2
	.
\end{align*}
The costs $C_l$ of the BIE solver consist of the two main components
\begin{enumerate}
	\item the cost $C^a_l$ of assembling the matrix of the linear system,
	\item the cost $C^s_l$ of solving the linear system.
\end{enumerate}

With the appropriate enumeration of the degrees of freedom, one can write the linear system in \change{\eqref{eq:lin_system_final}} as
\begin{align*}
	\begin{bmatrix}
    	\boldsymbol{\mathfrak{G}}_{11} & \boldsymbol{\mathfrak{G}}_{12}(\omega)
        \\
        \boldsymbol{\mathfrak{G}}_{21}(\omega) & \boldsymbol{\mathfrak{G}}_{22}(\omega) 
    \end{bmatrix}
    \begin{bmatrix}
    	\mu_1
        \\
        \mu_2
    \end{bmatrix}
    =
    \begin{bmatrix}
    	f_1
        \\
        f_2
    \end{bmatrix},
\end{align*}
where the block \change{$\boldsymbol{\mathfrak{G}}_{11}$} involves only points on the fixed boundary $\partial D_1$ and the remaining blocks depend on realizations of the boundary $\partial D_2(\omega)$.
Using the Schur complement, the above system can be reduced to the simpler one
\begin{align}\label{eq:lin_system_Schur}
	\Big( \boldsymbol{\mathfrak{G}}_{22}(\omega) - \boldsymbol{\mathfrak{G}}_{21}(\omega) \boldsymbol{\mathfrak{G}}_{11}^{-1} \boldsymbol{\mathfrak{G}}_{12}(\omega) \Big) \mu_2 &= f_2 - \boldsymbol{\mathfrak{G}}_{21}(\omega) \boldsymbol{\mathfrak{G}}_{11}^{-1} f_1,
    \\
    \mu_1 &= \boldsymbol{\mathfrak{G}}_{11}^{-1} f_1 - \boldsymbol{\mathfrak{G}}_{11}^{-1} \boldsymbol{\mathfrak{G}}_{12}(\omega) \mu_2.
\end{align}
Since $\boldsymbol{\mathfrak{G}}_{11}$ is fixed, it should be inverted only once using, e.g., LU decomposition.
Hence the cost of evaluating $\boldsymbol{\mathfrak{G}}_{11}^{-1}$ can be neglected.

Let $N_{1,l}$ and $N_{2,l}$ be the numbers of degrees of freedom corresponding to the fixed and random parts of the boundary $\partial D_1$ and $\partial D_2$ at the level $l \in [0,L]$.
Then the total number of degrees of freedom is $N_l=N_{1,l}+N_{2,l}$.
The construction of matrices $\boldsymbol{\mathfrak{G}}_{22}$, $\boldsymbol{\mathfrak{G}}_{21}$ and $\boldsymbol{\mathfrak{G}}_{12}$ requires $O(N_l^2-N_{1,l}^2)$ evaluation operations.
Evaluation of the Schur complement $\boldsymbol{\mathfrak{G}}_{22} - \boldsymbol{\mathfrak{G}}_{21} \boldsymbol{\mathfrak{G}}_{11}^{-1} \boldsymbol{\mathfrak{G}}_{12}$ requires $O( N_{1,l} N_{2,l} N_{l} + \change{N_{2,l}^2})$ arithmetic operations assuming standard matrix multiplication algorithm and that $\boldsymbol{\mathfrak{G}}_{11}^{-1}$ is given.
Evaluation of the right hand side requires additional $O(N_{1,l}N_{2,l} + N_{2,l})$ arithmetic operations since $\boldsymbol{\mathfrak{G}}_{11}^{-1} f_1$ is fixed. 
Hence the cost $C_l^a $ of assembling the matrix involves $O(N_{1,l} N_{2,l} N_{l} + N_{2,l}^2 + N_{1,l}N_{2,l} + N_{2,l})$ arithmetic operations and $O(N_l^2-N_{1,l}^2)$ evaluation operations.

For direct linear solvers, the cost of solving the resulting linear system is equal to $O(N_{2,l}^{\gamma})$ for some $\gamma \in (2,3]$; we assume that $\gamma = 3$.
Evaluation of $\mu_1$ requires only $O(N_{1,l}N_{2,l}+\change{N_{1,l}})$ arithmetic operations since $\boldsymbol{\mathfrak{G}}_{11}^{-1} f_1$ and $\boldsymbol{\mathfrak{G}}_{11}^{-1} \boldsymbol{\mathfrak{G}}_{12}(\omega)$ are available from the first equation.
Therefore, the overall cost of solving the linear system can be estimated as $C_l^s = O( N_{2,l}^3 + N_{1,l}N_{2,l} + N_{1,l} )$.

Omitting evaluation operations, one may conclude that the total cost at each level $l \in [0,L]$ behaves asymptotically as
\begin{align}\label{eq:CG_cost}
	{C_l} = C_l^a + C_l^{s} &= O \left( N_{2,l}^3 +  N_{1,l} N_{2,l} N_{l} + N_{2,l}^2 + N_{1,l}N_{2,l} + N_{l}\right)
	\simeq  h_l^{-3}
\end{align}
since $N_l \simeq |\partial D| h_l^{-1}$.

Taking into account \eqref{eq:CG_cost} and since $\rho=3$, $\beta=2\alpha\geq 4$, the $\epsilon$-cost of the MLMC estimator follows from Theorem \ref{the:I_MLMC_cost} as 
\begin{align}\label{eq:MLMC_cost}
	C_{\mathbb{ML}}^{FS}
	\simeq
	\change{\big( |\partial D_2|^3 + |\partial D_1| |\partial D_2| |\partial D| \big)} \epsilon^{-2}.
\end{align}

\subsubsection{Analytical Green's kernel.}
\label{sec:anGF_complexity}

For analytical Green's kernel, the cost $C_{l}^a$ of assembling the matrix involves only $O(N_{2,l}^2)$ evaluation operations and the cost of solving the linear system is $C_l^s = \change{O(N_{2,l}^{3})}$.
Omitting evaluation operations, the total asymptotical cost at each level $l \in [0,L]$ is just
\begin{align}\label{eq:cost_anal}
	C_l = C_l^a + C_l^{s} = O \left(  N_{2,l}^3 \right)
	\simeq |\partial D_2|^{3} h_l^{-3}.
\end{align}
Hence $\rho=3$ and the $\epsilon$-cost of the MLMC estimator follows from Theorem \ref{the:I_MLMC_cost} as
\begin{align}\label{eq:MLMC_cost_anGF}
	C_{\mathbb{ML}}^{GF}
	\simeq
	|\partial D_2|^3 \epsilon^{-2}.
\end{align}

\subsubsection{Approximate Green's kernel.}
\label{sec:numGF_complexity}

Both \eqref{eq:GF_num_reg_part} and \eqref{eq:lin_system_Schur} have  the Schur complement of $\boldsymbol{\mathfrak{G}}_{11}$ as their matrices.
Hence, the cost $C_{l}^a$ of assembling the matrix involves the same $O( N_{1,l} N_{2,l} N_{l} + N_{2,l}^2)$ arithmetic operations.
However, evaluation of the right hand side does not require additional operations.
The cost of solving the linear system is again $C_l^s = O(N_{2,l}^3)$, and the total cost at each level has the following asymptotical behavior 
\begin{align} \label{eq:cost_num}
	C_l = C_l^a + C_l^{s} &= O \left( N_{2,l}^{3} + N_{1,l} N_{2,l} N_{l} + N_{2,l}^2 \right)
	\simeq h_l^{-3}.
\end{align}

The overall asymptotical complexity of the MLMC estimator is thus given by
\begin{align}\label{eq:MLMC_cost_numGF}
	C_{\mathbb{ML}}^{GF}
	\simeq
	\big( |\partial D_2|^3 + |\partial D_1| |\partial D_2| |\partial D| \big) \epsilon^{-2}.
\end{align}

\textbf{Corollary.}
It is clear from \eqref{eq:MLMC_cost}, \eqref{eq:MLMC_cost_anGF} and \eqref{eq:MLMC_cost_numGF} that any choice of the kernel function results in the optimal asymptotical complexity of the MLMC estimator.
However, the direct comparison of \eqref{eq:CG_cost}, \eqref{eq:cost_anal} and \eqref{eq:cost_num} reveals that the classical method requires extra $O(N_{1,l} N_{2,l} N_{l} + N_{2,l}^2 + N_{1,l}N_{2,l} + N_{l})$ and $O(N_{1,l}N_{2,l} + N_{l})$ operations at each level when compared to the methods with analytical and numerical Green's kernels respectively.
In the next section, we provide the numerical assessment of the impact of the choice of the Green's kernel on the performance of the proposed algorithm.

\section{Numerical results}
\label{sec:numerical}

\paragraph{\bf{Example 1.}}

\begin{figure}[!t]
	\centering
    \includegraphics[width=0.4\textwidth]{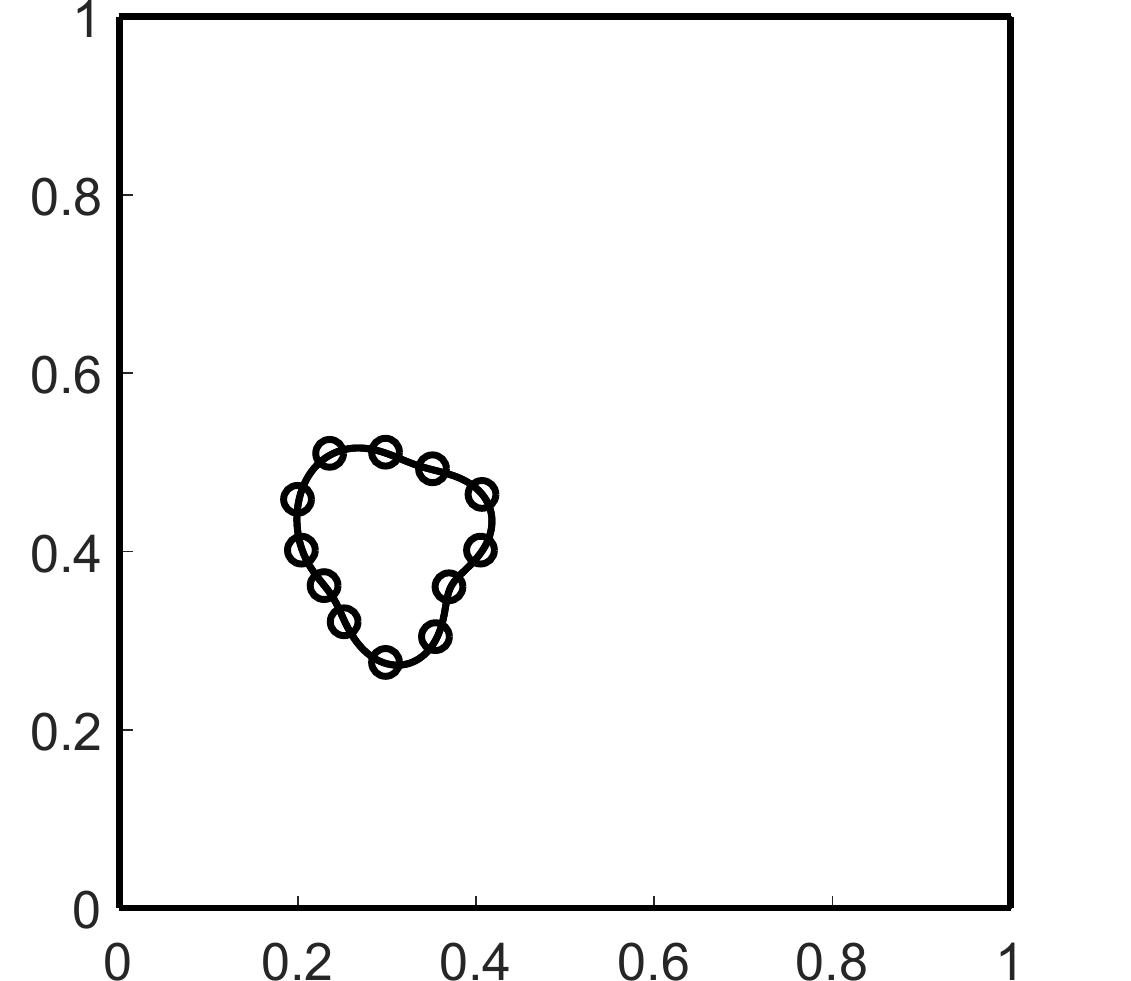}
    \qquad
    \includegraphics[width=0.4\textwidth]{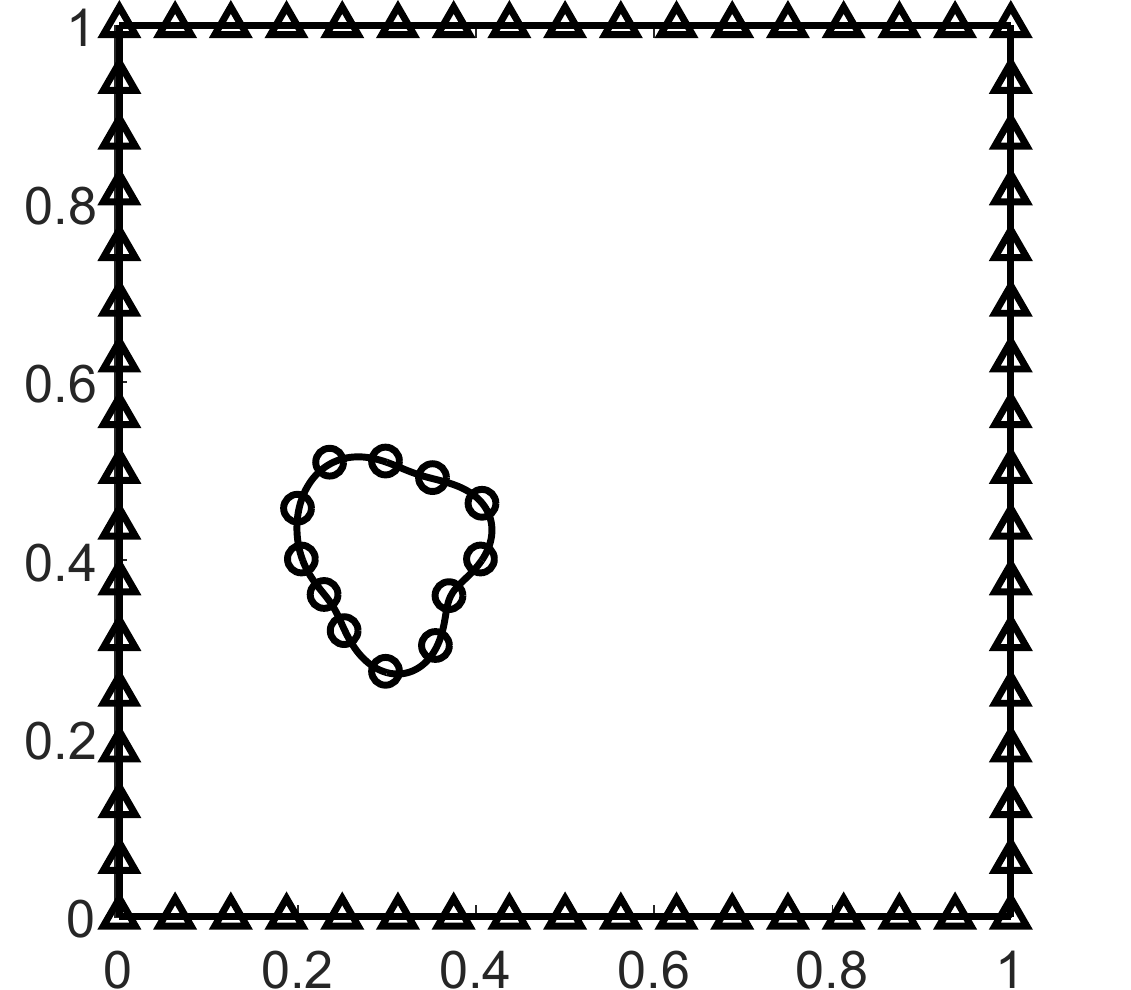}
	\caption[Discretization of the boundary.]{ Discretization of the boundary in Example 1 for BIE with analytical (left) and numerical (right) Green's kernels. }
	\label{fig:ex1_mesh_a}
\end{figure}

In the first example, we test the accuracy of the proposed discretization scheme for the fixed deterministic boundary.
Consider the problem
\begin{alignat*}{3}
	-\nabla^2 u(x) & = 0  &\qquad &\text{for } x \in D,
	\\ \nonumber
	u(x)           & = 0  &\qquad &\text{for } x \in \partial D_1,
	\\ \nonumber
	u(x)           & = G_1(x,\xi)  &\qquad &\text{for } x \in \partial D_2,
\end{alignat*}
where $D$ is the square domain with a single aperture. 
Figure \ref{fig:ex1_mesh_a} shows the discretization of the boundary for the case of analytical and numerical Green's kernels. 
For BIE with exact kernel, only the boundary $\partial D_2$ of the aperture has to be discretized while for BIE with approximate kernel, it is necessary to discretize both boundaries.
$G_1(x,\xi)$ is the Green's function \eqref{eq:spectral_GF_rect_summed} for the square bounded by $\partial D_1$ with the source located at the center of the aperture.
This choice of the boundary condition on $\partial D_2$ suggests $G_1(x,\xi)$ as the analytical solution of the above problem.

Convergence properties of the numerical scheme \eqref{eq:lin_system} with the quadrature nodes in \eqref{eq:colloc_points} are given in Tables \ref{tab:IV_Analit_GF_err} and \ref{tab:IV_Num_GF_err}.
The reference density of the potential $\mu$ was evaluated with the higher order scheme using the quadrature nodes in~(\ref{eq:colloc_points}').
The apparent rates of convergence are defined as $\alpha_h = \displaystyle{\frac{\log (e_{l}/e_{l-1})}{\log (h_l/h_{l-1})}}$.
It is seen that the errors have the order of convergence $\alpha=2$ in all norms as predicted by analysis.
Obviously, the errors in Table~\ref{tab:IV_Num_GF_err} have larger values but the difference is not large and the order is not reduced.

\begin{table}[!t]
\centering
\def\arraystretch{1.5}
\begin{tabular}{|c|cc|cc|cc|}
\hline
 $N_2$ & $\norm{\mu-\mu_h}{L^{\infty}}$ & rate & $\norm{u-u_h}{L^{\infty}}$ & rate & $\norm{u-u_h}{H^{1}}$ & rate \\ \hline 
   8 & 5.074e-1 &   -  & 3.323e-2 &   -  &  9.462e-2 &   -  \\ 
  20 & 7.008e-2 & 1.87 & 4.506e-3 & 2.36 &  2.992e-2 & 1.56 \\ 
  68 & 5.215e-3 & 2.06 & 1.820e-4 & 3.33 &  2.058e-3 & 2.64 \\
 260 & 3.533e-4 & 2.01 & 4.898e-6 & 2.34 &  1.312e-5 & 3.87 \\
1028 & 2.258e-5 & 2.00 & 3.133e-7 & 2.00 &  8.176e-7 & 2.00 \\
4100 & 1.428e-6 & 1.99 & 1.970e-8 & 2.00 &  5.140e-8 & 2.00 \\ \hline
\end{tabular}
\caption[ Convergence with analytical Green's kernel ]{ Convergence with analytical Green's kernel. }
\label{tab:IV_Analit_GF_err}
\end{table}

\begin{table}[!t]
\centering
\def\arraystretch{1.5}
\begin{tabular}{|cc|cc|cc|cc|}
\hline
  $N_1$ &  $N_2$ & $\norm{\mu-\mu_h}{L^{\infty}}$ & rate & $\norm{u-u_h}{L^{\infty}}$ & rate & $\norm{u-u_h}{H^{1}}$ & rate \\ \hline 
  48  &    8 & 5.071e-1 &   -  & 3.323e-2 &   -  &  1.048e-1 &   -  \\
 104  &   20 & 7.004e-2 & 1.87 & 4.507e-3 & 2.36 &  3.460e-2 & 1.53 \\
 328  &   68 & 5.212e-3 & 2.07 & 1.820e-4 & 3.33 &  2.462e-3 & 2.64 \\
1232  &  260 & 3.530e-4 & 2.01 & 4.901e-6 & 2.34 &  1.335e-5 & 3.98 \\
4848  & 1028 & 2.257e-5 & 2.00 & 3.135e-7 & 2.00 &  8.319e-7 & 1.99 \\
19336  & 4100 & 1.422e-6 & 1.99 & 1.971e-8 & 2.00 &  5.230e-8 & 2.00 \\ \hline
\end{tabular}
\caption[ Convergence with numerical Green's kernel ]{ Convergence with numerical Green's kernel. }
\label{tab:IV_Num_GF_err}
\end{table}

\begin{figure}[!t]
	\centering
    \includegraphics[width=0.25\linewidth]{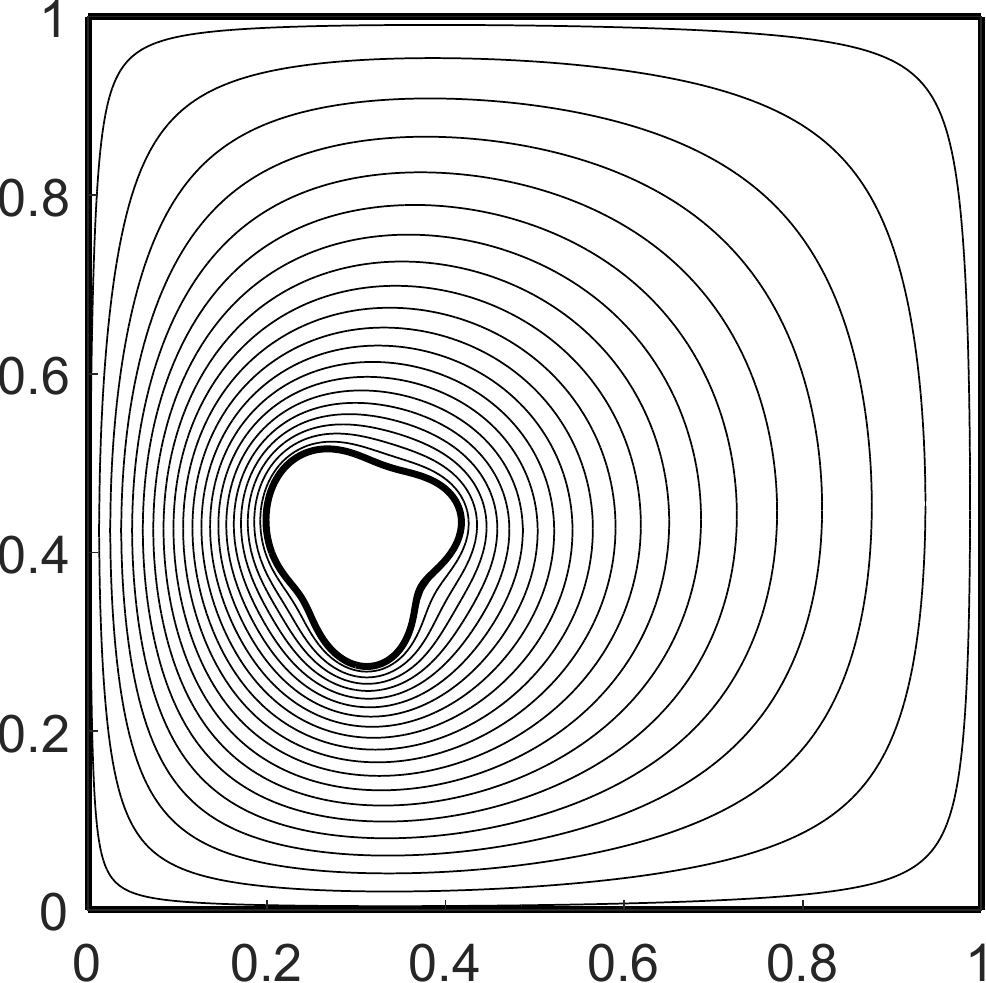}
    \;
    \includegraphics[width=0.32\textwidth]{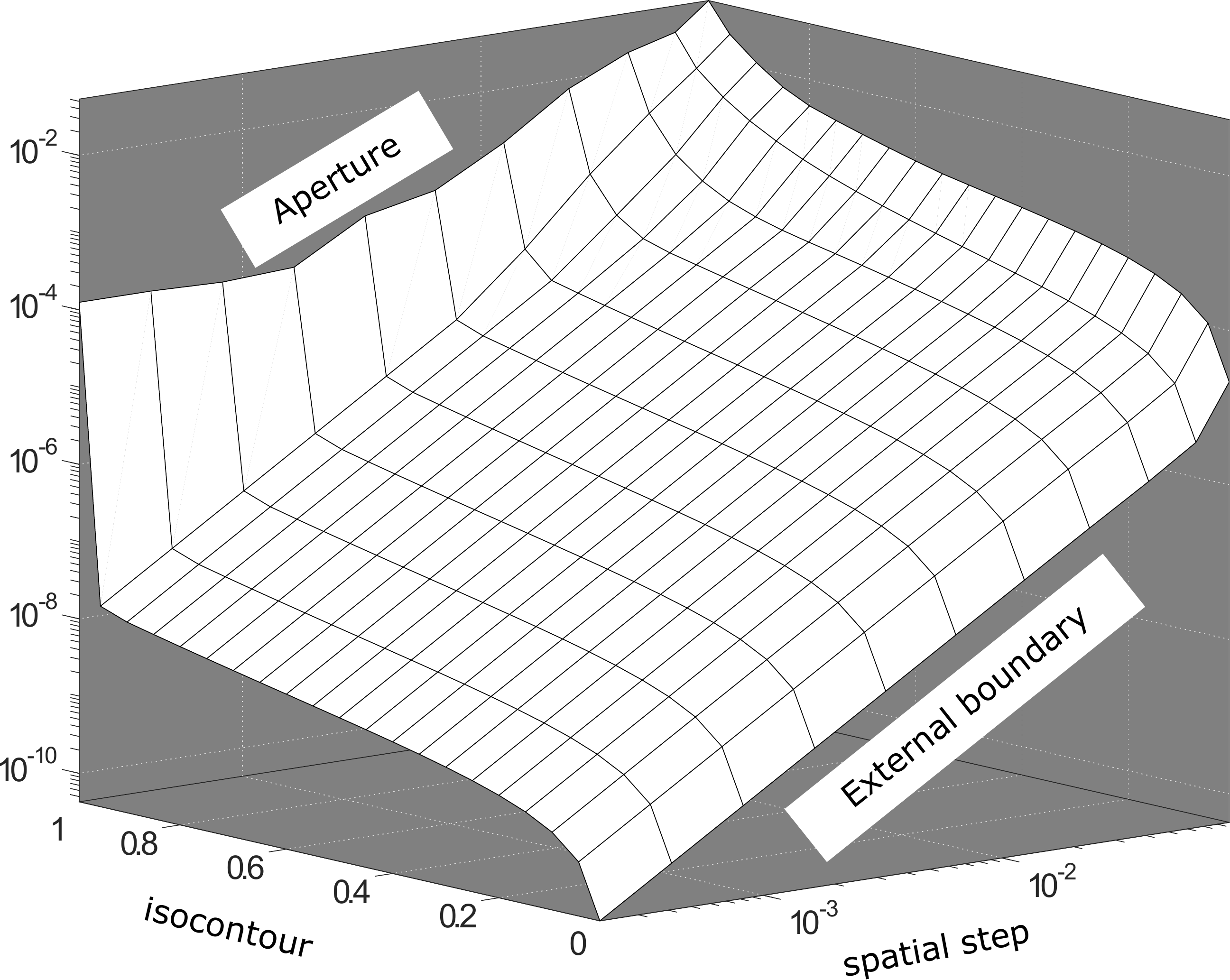}
    \;
    \includegraphics[width=0.32\textwidth]{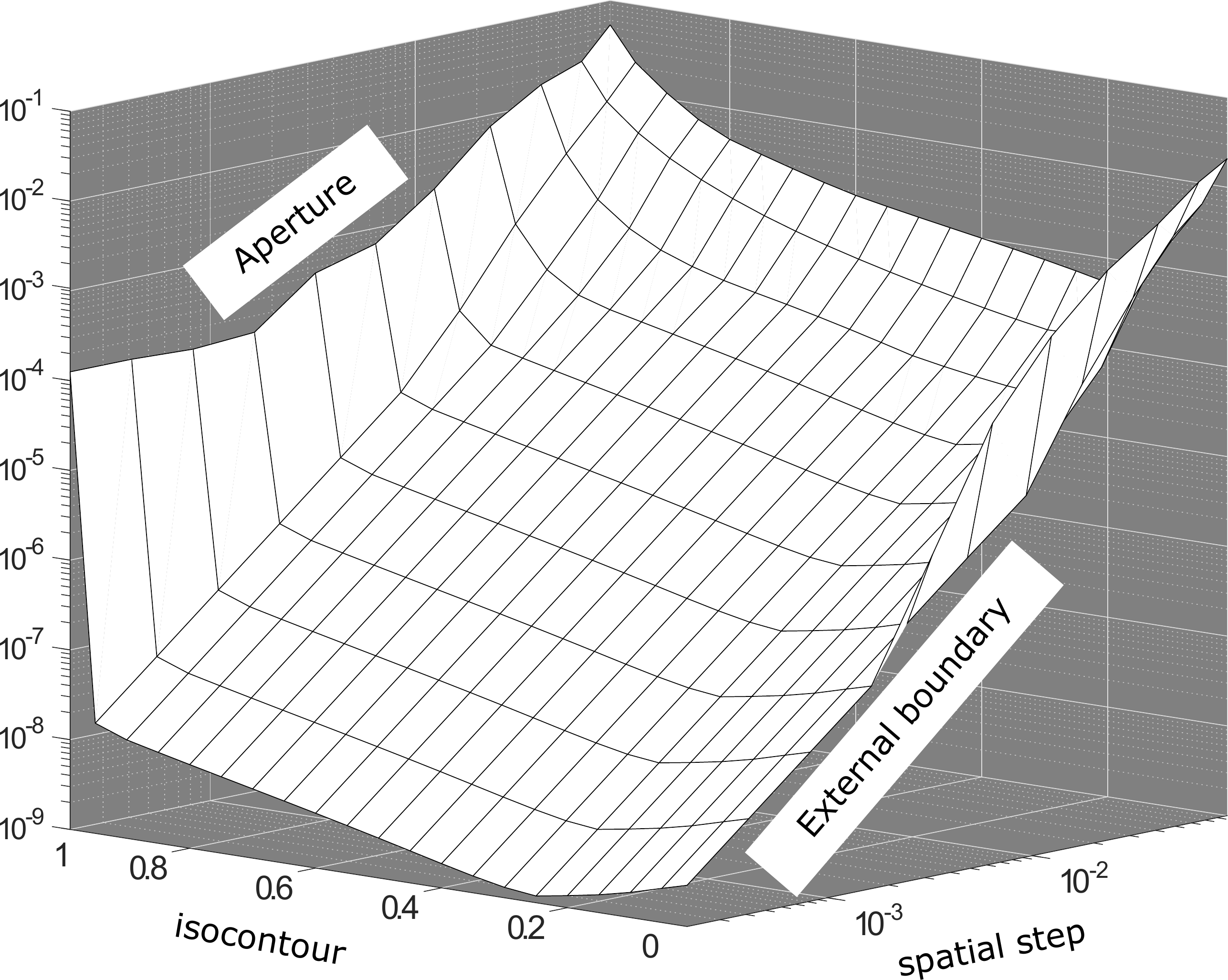}
	\caption[Maximum errors of the spatial approximation]{ Isocontours of the aperture (left) and $L_{\infty}$ errors along the isocontours for BIE with analytical (middle) and numerical (right) Green's kernels. }
	\label{fig:ex1_isol_conv}
\end{figure}

Results in Tables \ref{tab:IV_Analit_GF_err} and \ref{tab:IV_Num_GF_err} are also presented graphically in Figure \ref{fig:ex1_isol_conv} which illustrates the supremum norm of the error along the isocontours of the boundary.
Superiority of the analytical Green's kernel is obvious near the deterministic boundary $\partial D_1$ but both approaches show good results far from the boundaries.

\begin{figure}[!t]
	\centering
    \includegraphics[width=0.32\textwidth]{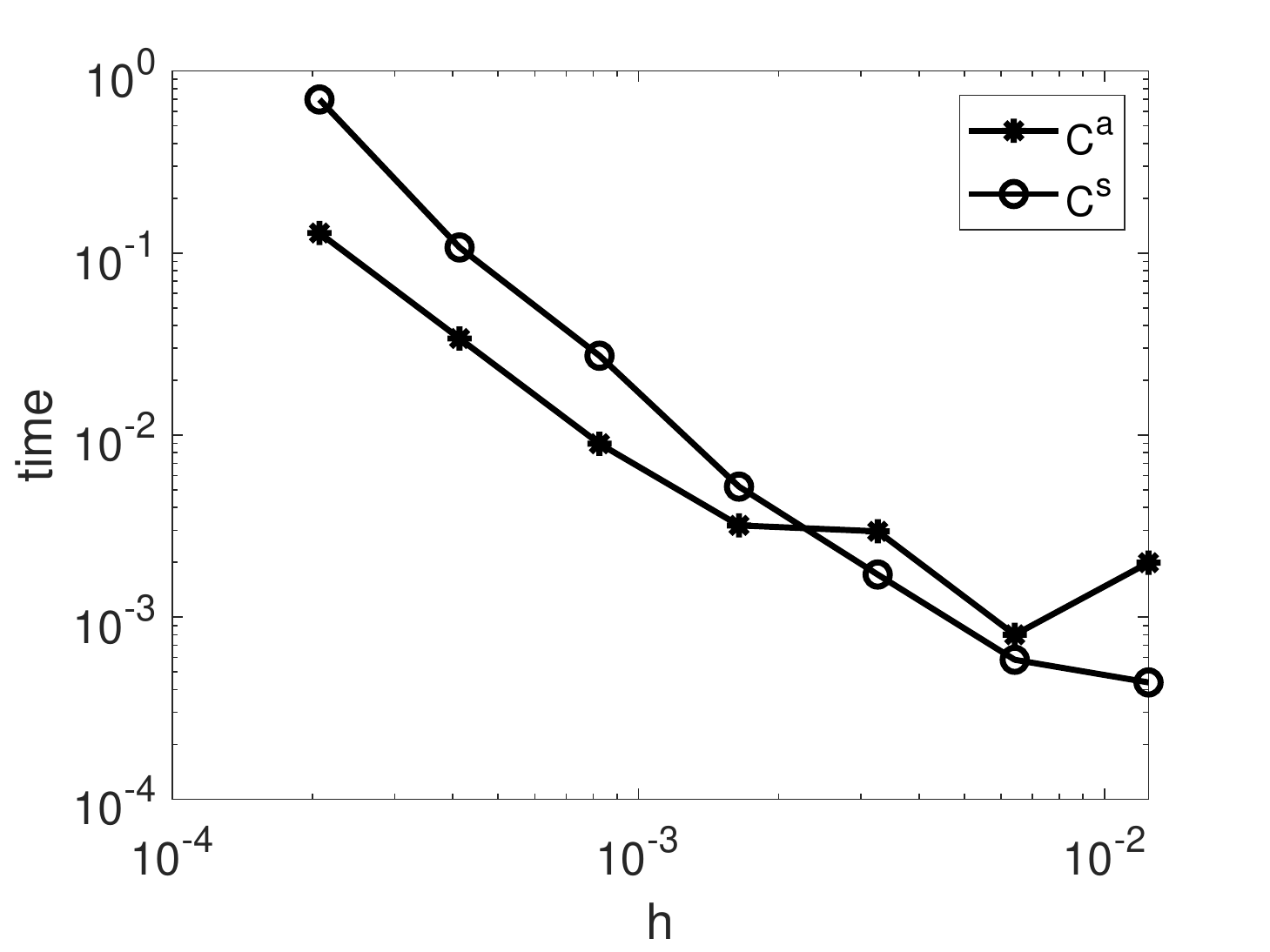}
    \includegraphics[width=0.32\textwidth]{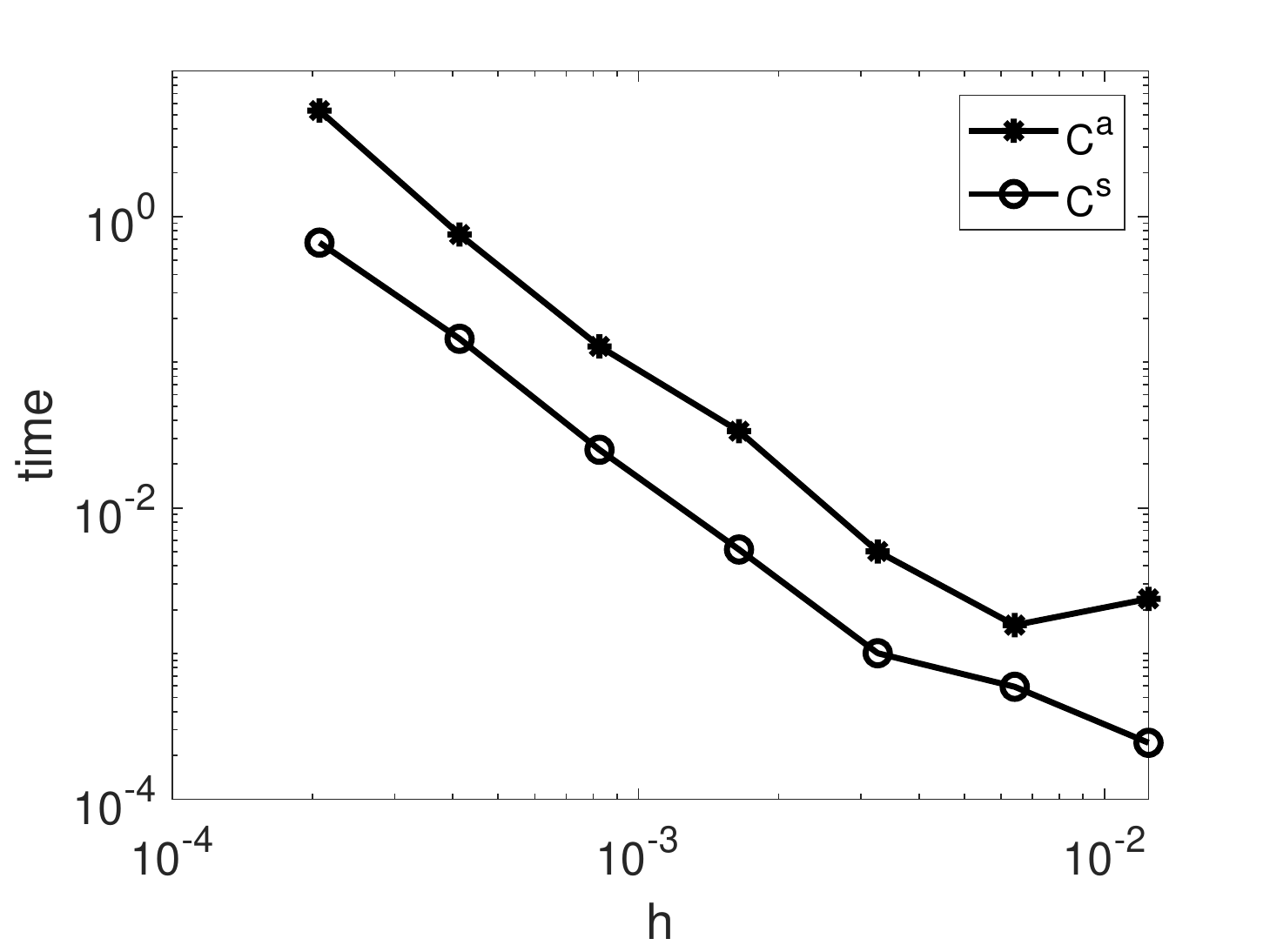}
    \includegraphics[width=0.32\textwidth]{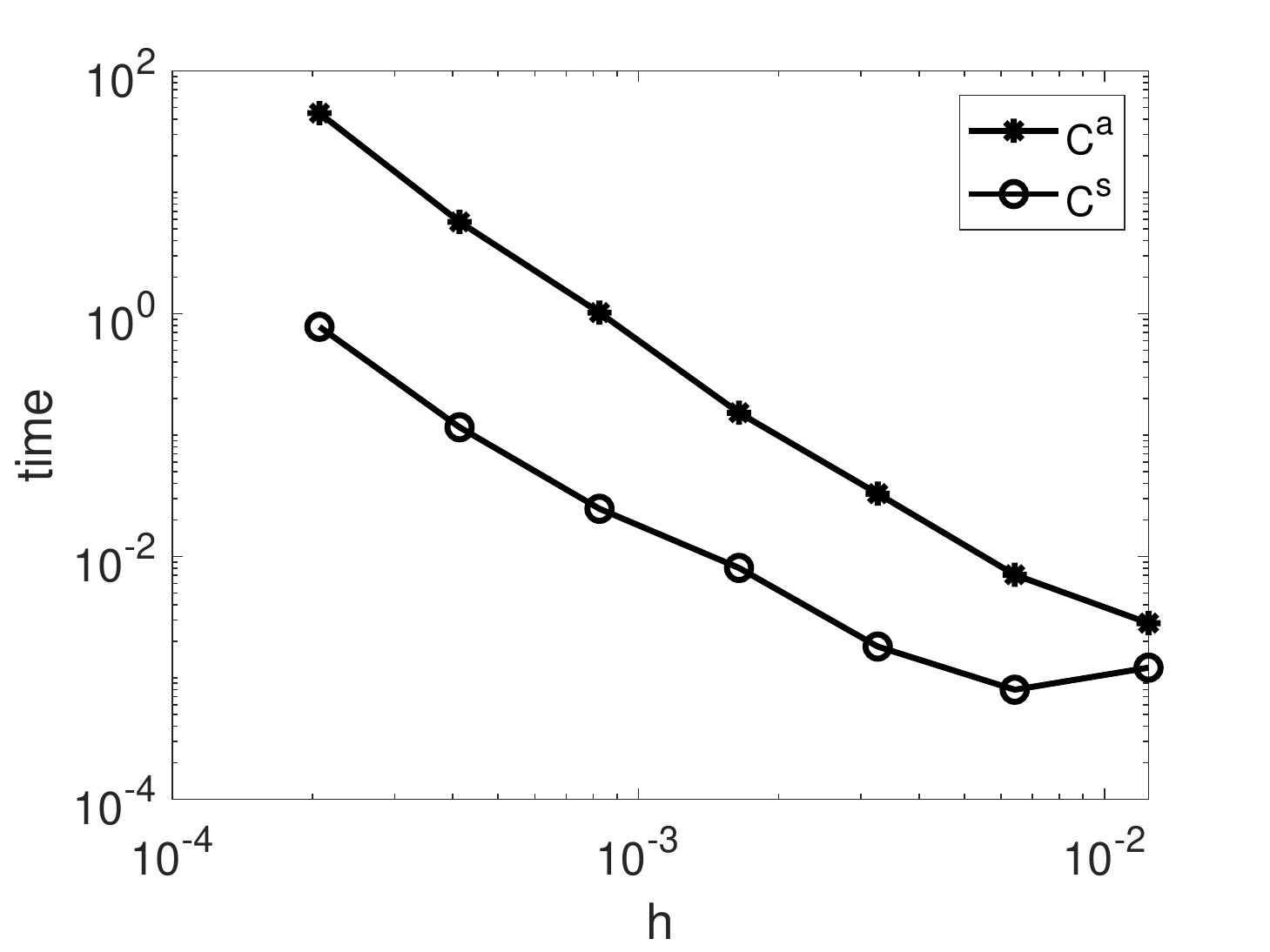}
	\caption[The costs $C^a$ and $C^s$ of assembling and solving the linear system]{ The costs $C^a$ and $C^s$ of assembling and solving the linear system for the schemes with analytical Green's kernel (left), numerical Green's kernel (middle) and fundamental solution (right).}
	\label{fig:ex1_cost_per_step}
\end{figure}

\begin{figure}[!t]
	\centering
    \includegraphics[width=0.4\textwidth]{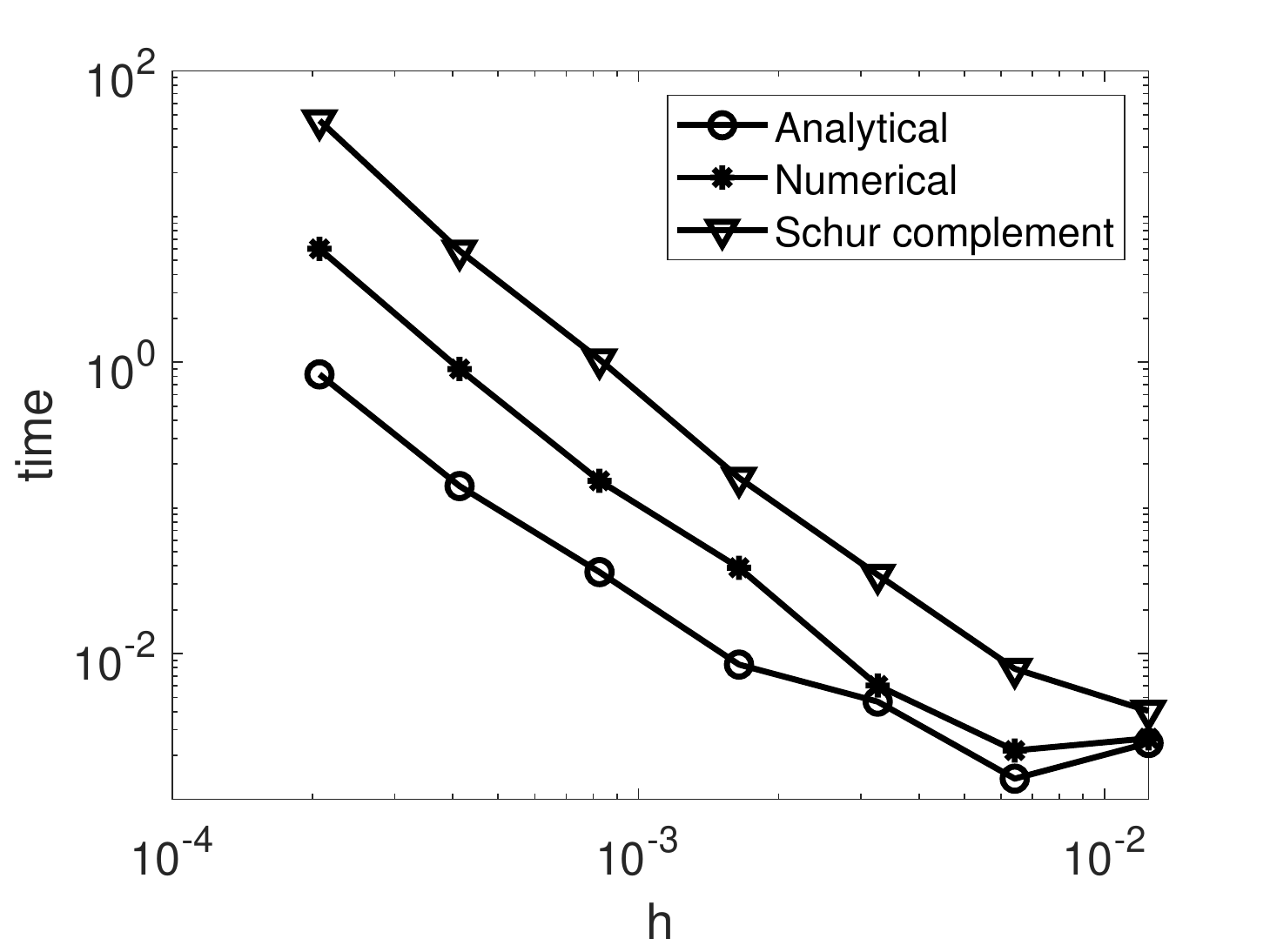}
    \quad
    \includegraphics[width=0.4\textwidth]{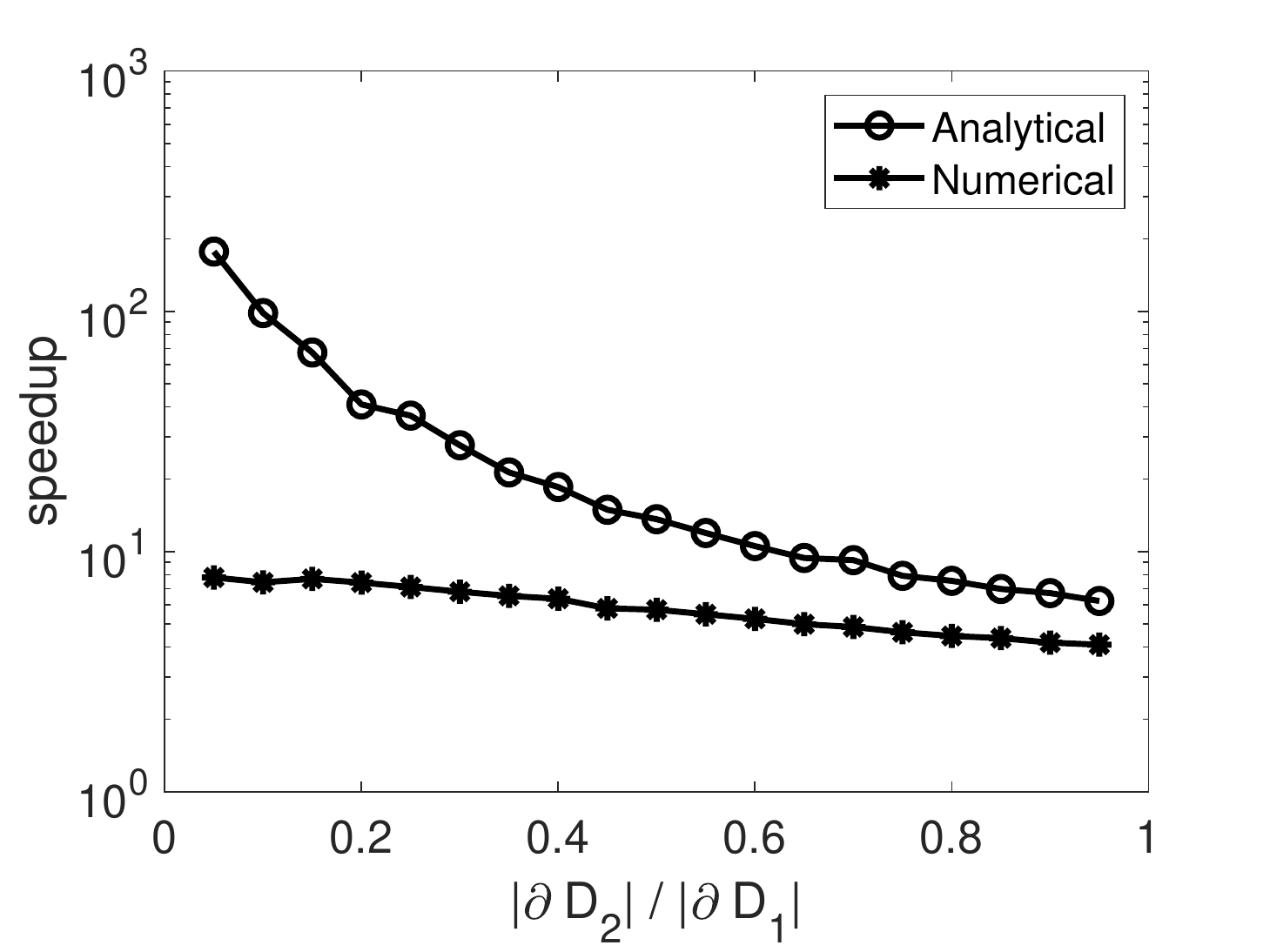}
	\caption{ Empirical costs $C_l$ in \eqref{eq:CG_cost}, \eqref{eq:cost_anal}, \eqref{eq:cost_num} (left) and corresponding computational speedups of the schemes with analytical and numerical kernels (right).}
	\label{fig:ex1_speedup}
\end{figure}

Figure \ref{fig:ex1_cost_per_step} illustrates the costs $C^a$ of assembling the matrix and $C^s$ of solving the resulting system.
As predicted, the cost $C^s$ is asymptotically dominant for BIE with analytical Green's kernel while the cost $C^a$ dominates in the case of numerical Green's function.
Figure \ref{fig:ex1_speedup} depicts the overall empirical computational costs of different methods along with the speedup of the proposed schemes over the standard BIE method with fundamental solution as kernel.
It is seen that methods with both analytical and numerical Green's kernels are superior to standard approach and perform exceptionally well for small ratios of $\frac{|\partial D_2|}{|\partial D_1|}$.

\paragraph{\bf{Example 2.}}
For the second example, consider the problem
\begin{alignat*}{3}
	-\nabla^2 u(x) & = f(x)  &\qquad &\text{for } x \in D,
	\\ \nonumber
	u(x)           & = 0  &\qquad &\text{for } x \in \partial D_1,
	\\ \nonumber
	\frac{\partial u(x)}{\partial n}           & = 0  &\qquad &\text{for } x \in \partial D_2,
\end{alignat*}
where $D$ and $\partial D_1$ are the same as in Example 1 and $\partial D_2$ is the aperture with the following parametrization
\begin{align*}
	\big(x(t,\omega),y(t,\omega)\big) = \big(x_c(\omega),y_c(\omega)\big) + R(t,\omega) \big(\cos(t),\sin(t)\big), \quad  t \in [0,2 \pi).
\end{align*}
The radius of the aperture is defined as
\begin{align*}
	R(t,\omega) = \overline{R}(t) + \sigma_r \sum_{n=1}^s \Big( a_n(\omega) \cos(2\pi n t) + b_n(\omega) \sin(2\pi n t) \Big)
\end{align*}
with the mean radius $\overline{R}(t)$ and the random coefficients $a_n(\omega) = \mathcal{U}(-\sqrt{3},\sqrt{3})$, $b_n(\omega) = \mathcal{U}(-\sqrt{3},\sqrt{3})$.
Coefficient $\sigma_r$ controls intensity of the random perturbation.
Coordinates of the center are also random variables
\begin{align*}
	x_c(\omega) &= \overline{x}_c + \sigma_x \mathcal{U}(-1,1),
	\\ \nonumber
	y_c(\omega) &= \overline{y}_c + \sigma_y \mathcal{U}(-1,1),
\end{align*}
where $(\overline{x}_c,\overline{y}_c)$ is the mean location of the center and the coefficients $\sigma_x$, $\sigma_y$ control deviation from the mean.
For this particular example, we set $(\overline{x}_c,\overline{y}_c) = (0.3,0.4)$, $\sigma_x=\sigma_y=0.05$, $\overline{R}=0.15$, $\sigma_r=0.01$ and $s=10$.

\begin{figure}[!t]
	\centering
    \includegraphics[width=0.28\textwidth]{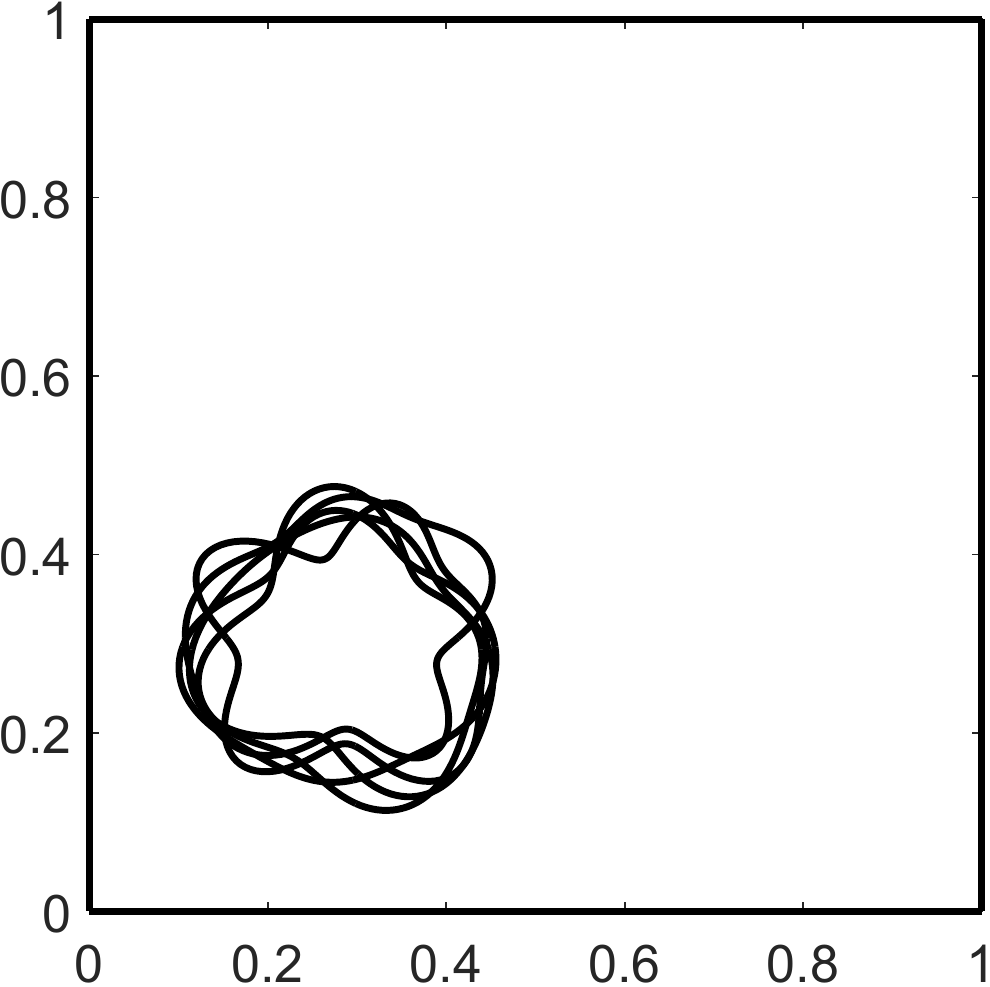}
    \;
    \includegraphics[width=0.28\textwidth]{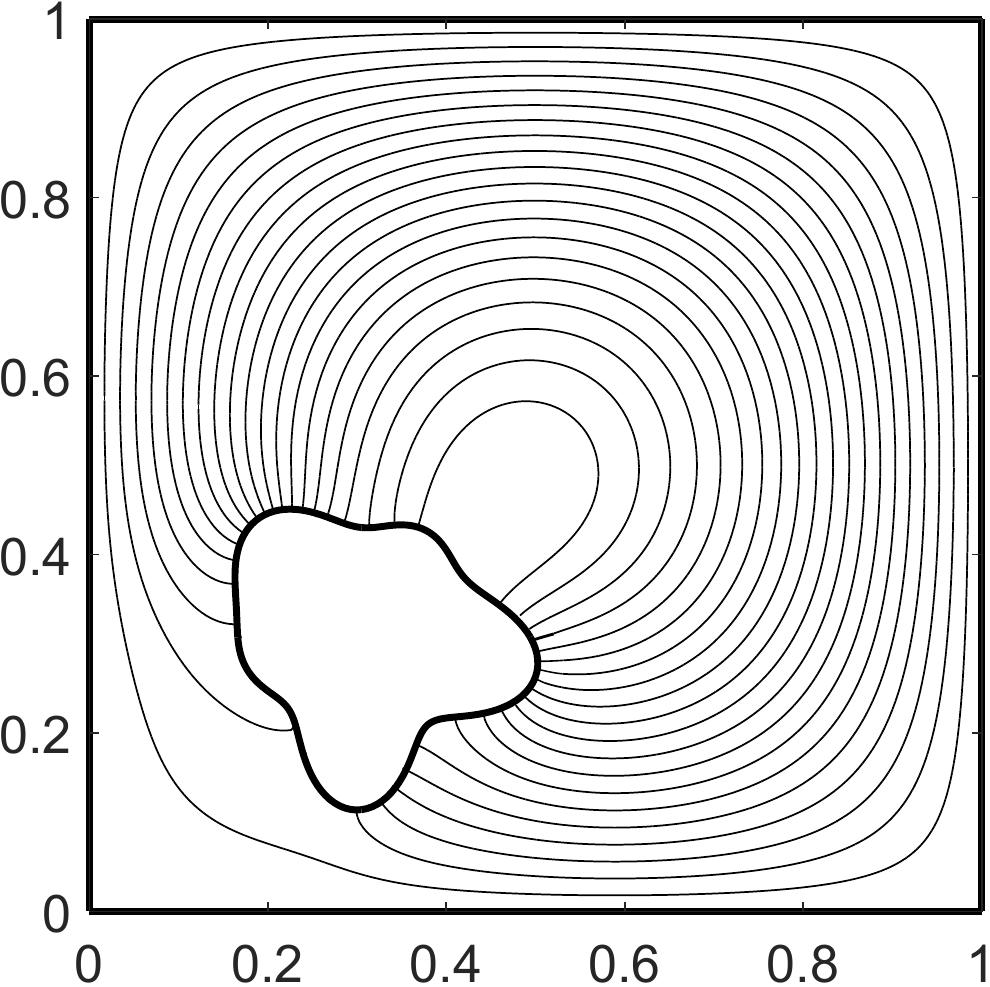}
    \;
    \includegraphics[width=0.28\textwidth]{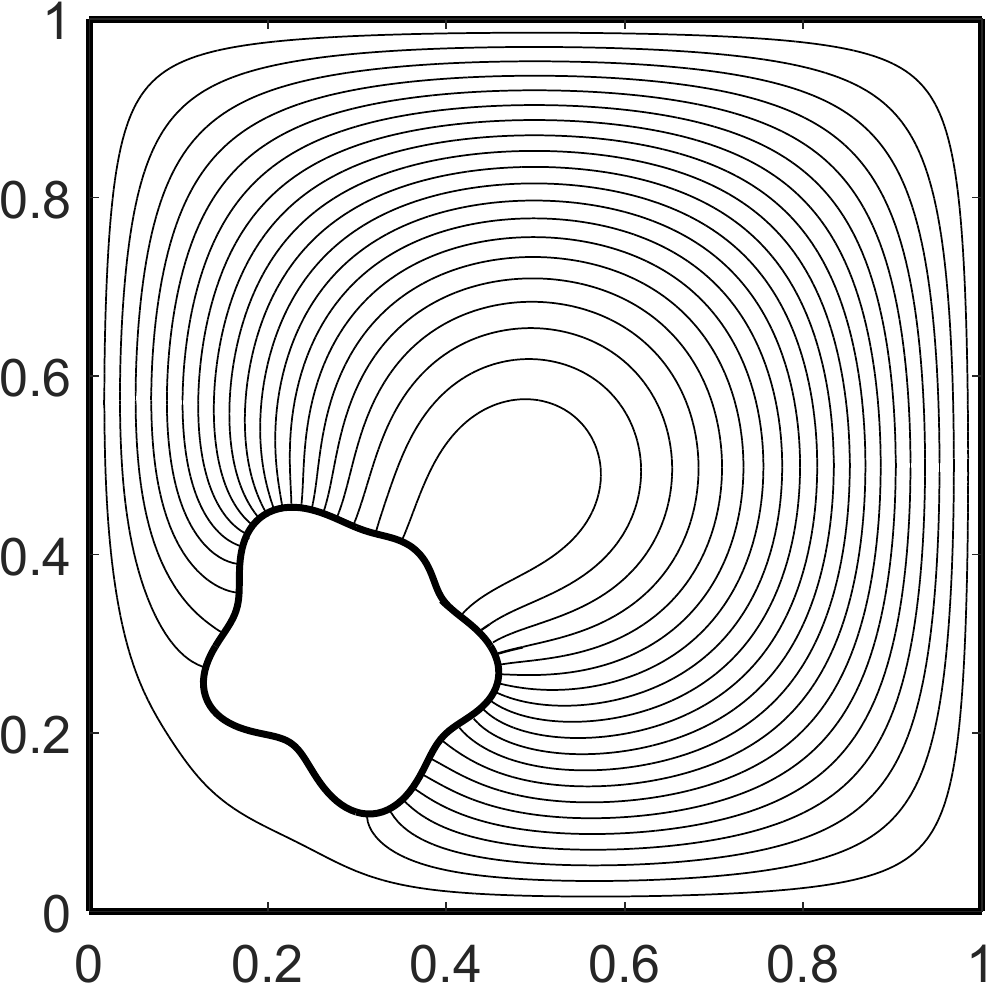}
	\caption{ Realizations of the random geometry (left) and isolines of the solutions corresponding to two different realizations.}
	\label{fig:rand_geom}
\end{figure}

The forcing term $f(x)$ is chosen such that the analytical solution in the deterministic domain $D_1$ without the aperture is given by
\begin{align*}
    u_1(x) = 
    100 \sum_{i=n}^2 \sum_{m=1}^2 \frac{ \sin(n \pi /2)^2 \sin(m \pi /2)^2 }{ nm\pi^4 (n^2+m^2) }
    \sin(n\pi x) \sin(m\pi y).
\end{align*}
Figure~\ref{fig:rand_geom} illustrates five different realizations of the random geometry and isolines of the solutions corresponding to two particular realizations.

\begin{figure}[!t]
	\centering
    \includegraphics[width=0.35\textwidth]{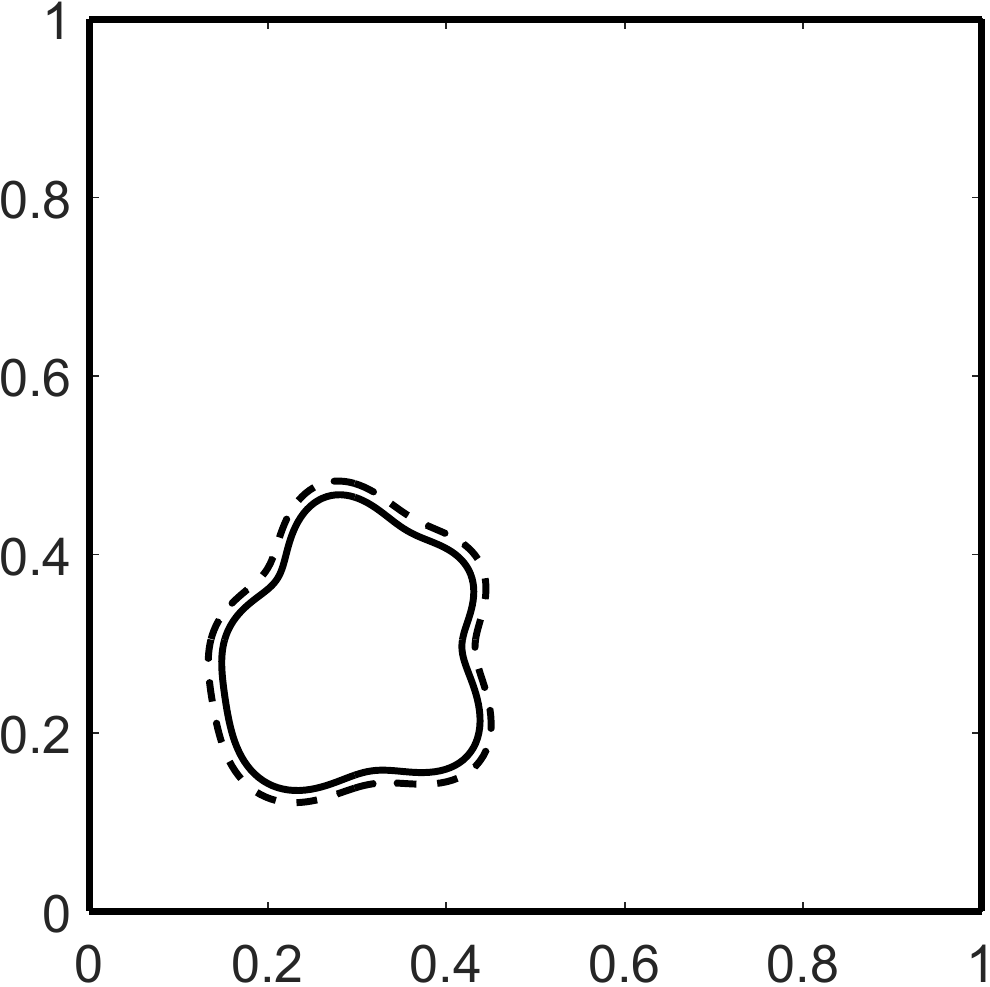}
    \qquad
    \includegraphics[width=0.45\textwidth]{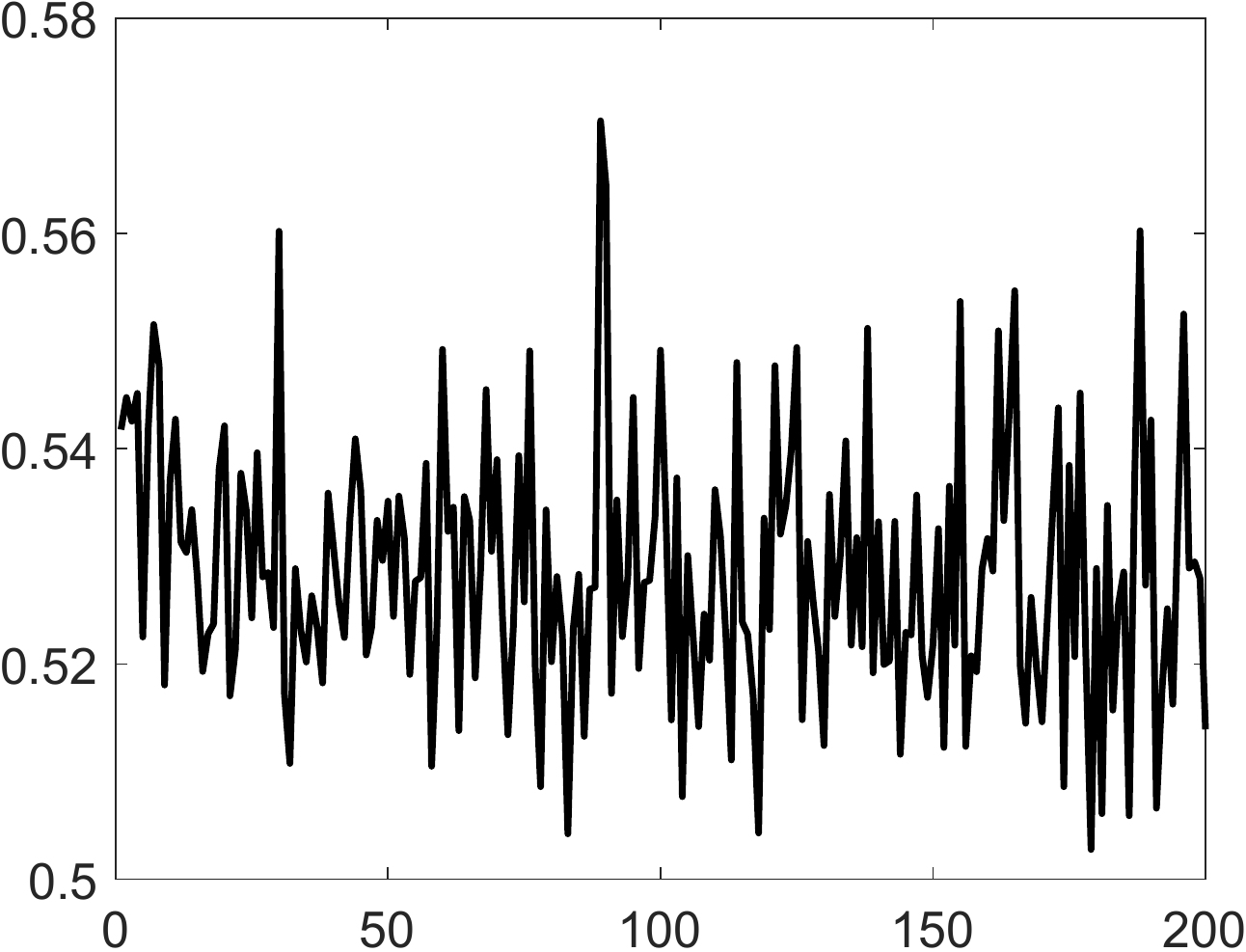}
	\caption{ Realization of the geometry and the contour $C$ (left) and 200 realizations of the functional (right) in \eqref{eq:functional}.}
	\label{fig:functional}
\end{figure}

Consider the functional 
\begin{align}\label{eq:functional}
    F[u]:=\sup_{x \in C} u(x),
\end{align}
where $C$ is the contour parallel to the aperture with the offset $d=0.01$.
A single realization of the geometry and the corresponding contour $C$ is depicted on Figure~\ref{fig:functional} along with 200 realizations of the value of the functional.

\begin{table}[t]
\centering
\begin{tabular}{|c|c|c|c|}
\hline
         & Analytical GF & Numerical GF & Schur complement \\ \hline
$\alpha$ & $1.88$        & $2.02$       & $1.89$           \\
$\beta$  & $4.33$        & $4.31$       & $4.42$           \\
$\rho$ & $1.58$        & $1.97$       & $2.52$           \\ \hline
\end{tabular}
\caption{Empirical rates $\alpha$ of the weak convergence, $\beta$ of the decay of the variance and $\rho$ of the growth of the cost for boundary integral equations with different kernels. }
\label{tab:mlmc_rates}
\end{table}

\begin{figure}[!t]
	\centering
    \includegraphics[width=0.45\textwidth]{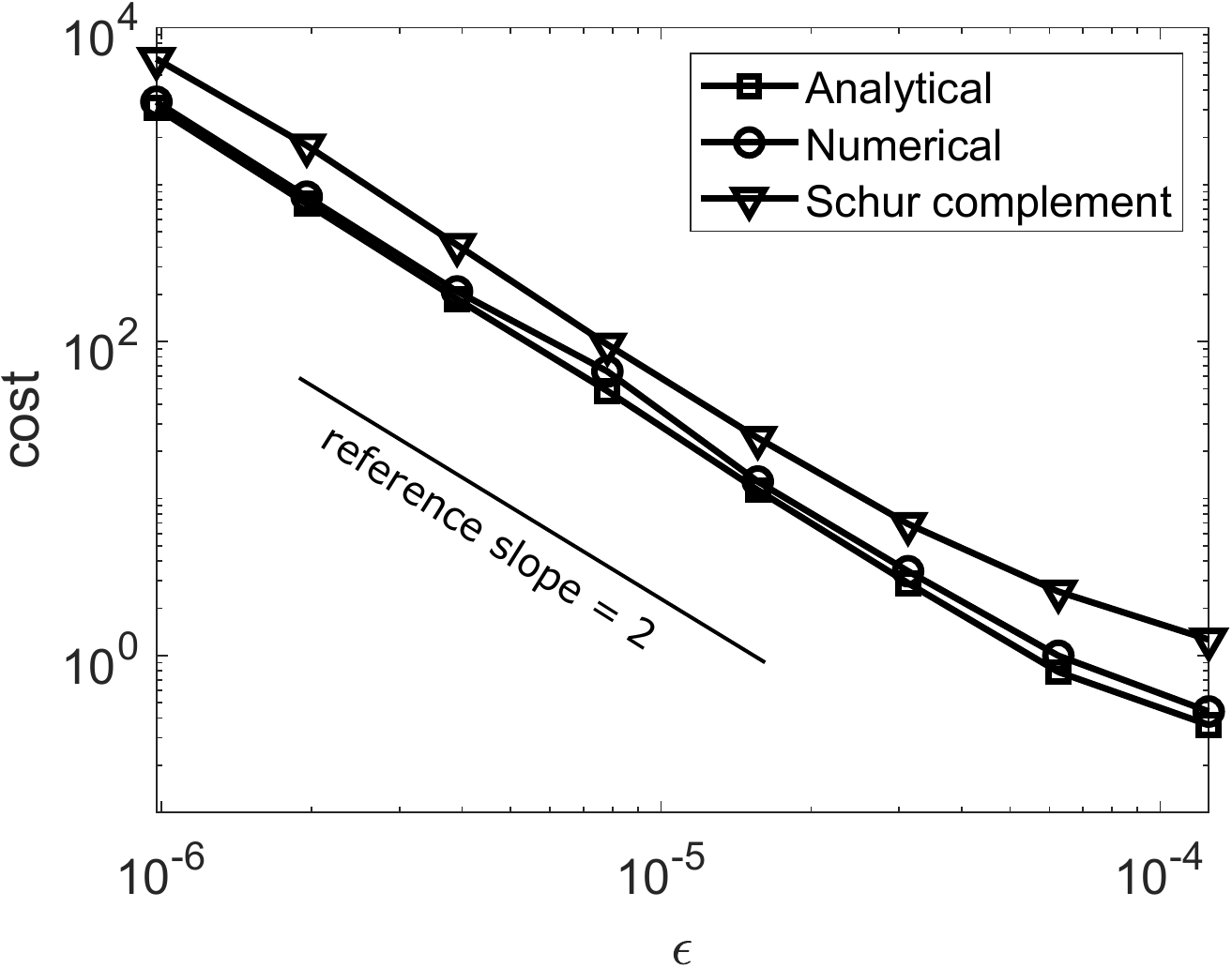}
	\caption{ $\epsilon$-cost of the MLMC method for the boundary integral equations with different kernels.}
	\label{fig:mlmc_cost}
\end{figure}

The observed rates $\alpha$, $\beta$, $\rho$ in Theorem \ref{the:I_MLMC_cost} for \change{the proposed numerical scheme} are given in Table \ref{tab:mlmc_rates}.
It is seen that, due to the second order of convergence of the spatial approximation, the variance decay is much faster than the growth of the cost for all types of kernels.
Therefore, the theoretical $\epsilon$-complexity of the MLMC estimator is proportional to $\epsilon^{-2}$ and Figure \ref{fig:mlmc_cost} confirms this prediction. 
One can also see that the proposed schemes with both analytical and numerical Green's kernels perform better than the standard method in section~\ref{sec:Schur_complexity} as expected.

\section{Conclusion}

In this paper, we considered the problem of approximating solutions to elliptic PDEs in domains comprised of deterministic and random boundaries.
We used the numerical scheme based on boundary integral representation of solutions to such PDEs and proposed to use Green's functions as the kernels of the corresponding potentials.
The proposed numerical scheme can be applied to problems in arbitrary domains \change{with random apertures} and does not require the explicit knowledge of analytical Green's functions which can be considered as a major contribution of this work.
We showed that when a large number of repetitive solutions is required, as is the case of Monte Carlo simulations, the proposed scheme can lead to significant reduction of computational complexity compared to standard BIE techniques.

\section*{Data availability}

The authors declare that all data supporting the findings of this study are available within the article.

\section*{Acknowledgements}

This material is based upon work supported in part by: the U.S. Department of Energy, Office of Science, Early Career Research Program under 
award number ERKJ314; U.S. Department of Energy, Office of Advanced Scientific Computing Research under award numbers ERKJ331 and ERKJ345; 
the National Science Foundation, Division of Mathematical Sciences, Computational Mathematics program under contract number DMS1620280;
and by the Laboratory Directed Research and Development program at the Oak Ridge National Laboratory, which is operated by UT-Battelle, LLC., 
for the U.S. Department of Energy under contract DE-AC05-00OR22725.


\bibliographystyle{spmpsci}      


\bibliography{main}

\end{document}